\documentclass[11pt,a4paper]{article}

\usepackage[utf8]{inputenc}
\usepackage[margin=2.5cm]{geometry}
\usepackage{amsmath,amssymb,amsthm}
\usepackage{booktabs}
\usepackage{array}
\usepackage{enumitem}
\usepackage[colorlinks=true,linkcolor=blue,citecolor=blue,urlcolor=blue]{hyperref}

\theoremstyle{plain}
\newtheorem{theorem}{Theorem}[section]
\newtheorem{lemma}[theorem]{Lemma}
\newtheorem{proposition}[theorem]{Proposition}
\newtheorem{corollary}[theorem]{Corollary}

\theoremstyle{definition}
\newtheorem{definition}[theorem]{Definition}
\theoremstyle{remark}
\newtheorem{remark}[theorem]{Remark}

\newcommand{\NF}{N_F}
\newcommand{\rNF}{\rho_{N_F}}
\newcommand{\rUps}{\rho_{\Upsilon}}
\newcommand{\dmax}{d_{\max}}
\newcommand{\dmin}{d_{\min}}
\newcommand{\demax}{\delta_{\max}}
\newcommand{\demin}{\delta_{\min}}
\newcommand{\E}{\mathcal{E}}
\newcommand{\one}{\mathbf{1}}
\newcommand{\bxi}{\boldsymbol{\xi}}
\newcommand{\bpi}{\boldsymbol{\pi}}
\newcommand{\Sp}[1]{S_{#1}^{+}}
\DeclareMathOperator{\tr}{tr}
\DeclareMathOperator{\rank}{rank}
\DeclareMathOperator{\diag}{diag}

\title{Neighbourhood second Zagreb matrix of graphs: spectral radius bounds and extremal graphs}

\author{%
  Subhendu Chakraborty\textsuperscript{1,2,}\thanks{\texttt{subhendu.chakraborty.math@gmail.com}},~
  Priti Prasanna Mondal\textsuperscript{3,}\thanks{\texttt{pritiprasanna1992@gmail.com}}\\[2ex]
  {\small \textsuperscript{1}Institute for Advancing Intelligence (IAI), TCG CREST, Kolkata-700091, India}\\
  {\small \textsuperscript{2}Academy of Scientific and Innovative Research (AcSIR), Ghaziabad-201002, India}\\
  {\small \textsuperscript{3}Indian Institute of Technology Hyderabad, Kandi, Sangareddy-502284, Telangana, India}%
}

\date{}
\begin{document}
\maketitle

\begin{abstract}
Let $G$ be a connected graph on at least two vertices, and let $\delta(i)$ be the sum of the degrees of the neighbours of a vertex $i$. The \emph{neighbourhood second Zagreb matrix} $\NF(G)$ has $(i,j)$-entry $\delta(i)\delta(j)$ if $i\sim j$ and $0$ otherwise; its spectral radius is denoted by
$\rNF(G)$. We prove that $\rNF(G)\ge\rho(G)^{5}$, where $\rho(G)$ is the adjacency spectral
radius, with equality if and only if $G$ is regular or bipartite semiregular, and we use this to
show that the same graphs are exactly the extremal graphs of a Berman--Zhang--type upper bound. We sharpen the estimates $\demin^{2}\rho(G)\le\rNF(G)\le\demax^{2}\rho(G)$ by means of the extreme edge weights of $\NF(G)$, give sharp upper bounds in terms of the maximum degree $\dmax$, the maximum neighbourhood degree sum $\demax$, the size $m$ and the clique number $\omega(G)$, and determine the extremal graphs for bipartite graphs and trees. For every connected unicyclic graph $G$ of order $n\ge3$ we prove
$32=\rNF(C_n)\le\rNF(G)\le\rNF(\Sp n)$, where $\Sp n$ is the triangle with $n-3$ pendant edges at
one vertex, with equality only for $C_n$ and only for $\Sp n$, respectively. Finally, we give bounds in terms of the first Zagreb index and for the
neighbourhood second Zagreb energy; for lower energy bounds, the extremal graphs are described completely.
\end{abstract}


\section{Introduction}

All graphs considered here are finite, simple, and undirected. Let $G$ be a graph with vertex set
$V(G)=\{1,2,\dots,n\}$ and edge set $E(G)$, $|E(G)|=m$. If two vertices $i$ and $j$ are adjacent,
we write $i\sim j$ and denote the corresponding edge by $\{i,j\}$. The \emph{adjacency matrix} of
$G$ is $A(G)=[a_{ij}]$, where $a_{ij}=1$ if $i\sim j$ and $a_{ij}=0$ otherwise. Since $A(G)$ is
real and symmetric, its eigenvalues $\mu_1\ge\mu_2\ge\dots\ge\mu_n$ are real; the largest one,
$\mu_1(G)$, is the \emph{spectral radius} of $G$ and is denoted by $\rho(G)$.

For $i\in V(G)$ let $N_G(i)$ be the neighbourhood of $i$, $N_G[i]=N_G(i)\cup\{i\}$ its closed
neighbourhood and $d(i)=|N_G(i)|$ its degree. Following the standard terminology of chemical graph
theory, the \emph{neighbourhood degree sum} of $i$ is
\[
\delta(i)=\sum_{j\sim i}d(j).
\]
We write $\demax=\demax(G)=\max\{\delta(i):i\in V(G)\}$ and
$\demin=\demin(G)=\min\{\delta(i):i\in V(G)\}$, and similarly $\dmax$, $\dmin$ for the maximum and
minimum degree. Clearly
\begin{equation}\label{eq:delta-range}
\dmin^{2}\le\demin\le\demax\le\dmax^{2}.
\end{equation}

Degree-based graph invariants such as the first and second Zagreb indices \cite{GT72, GRTW75} and the
forgotten index \cite{FG15} have been studied intensively, and in recent years their neighbourhood
versions, in which $d(i)$ is replaced by $\delta(i)$, have attracted considerable attention
\cite{MDP19,MDP21}. Ghorbani and Hosseinzadeh \cite{GH10} introduced the neighbourhood first and
second Zagreb indices
\[
M_1'(G)=\sum_{\{i,j\}\in E(G)}\bigl(\delta(i)+\delta(j)\bigr),\qquad
M_2'(G)=\sum_{\{i,j\}\in E(G)}\delta(i)\delta(j).
\]
The index $M_2'$ also appears as  $M_2^{*}$ \cite{MDP19} and as $NM_2$ in the literature. Note that $M_{N}(G)=\sum_{v}\delta(v)^2$ of \cite{MDP19, MDP21} differs from $M_{N}^{'}(G)=\sum_{v}d(v)\delta(v)$. We use the index $M_2'$ in Theorem~\ref{thm:KM}.

In parallel, matrices whose nonzero entries are degree-based edge weights have been introduced in order to study such invariants by spectral methods. Given a symmetric bivariate function $F$, the
associated weighted adjacency matrix $A_F(G)$ has $(i,j)$-th entry $F(d(i),d(j))$ if $i\sim j$ and
$0$ otherwise \cite{DGMMF18}; the extremal theory of $\rho(A_F)$ over trees and unicyclic graphs has
been developed by Li and Wang \cite{LW21} and by Zheng, Guan and Jin \cite{ZGJ23} under
monotonicity and convexity hypotheses on $F$.
Jafari Rad, Jahanbani and Gutman \cite{JJG18} considered in particular the first and second Zagreb
matrices, whose $(i,j)$-th entries are $d(i)+d(j)$ and $d(i)d(j)$ on the edges of $G$; Mondal et al. \cite{MBDP22} replaced $d(i)$ by $\delta(i)$ in the additive case, studied the energy of the
resulting neighbourhood first Zagreb matrix $\Upsilon(G)=[\delta(i)+\delta(j)]$ and established its value as a molecular
descriptor. The multiplicative counterpart, which is the matrix considered in the present paper,
was introduced by Barik and Verma \cite{BV25}.

\begin{definition}[\cite{BV25}]
The \emph{neighbourhood second Zagreb matrix} of a graph $G$, denoted by $\NF(G)$ (or simply
$\NF$), is the $n\times n$ matrix $\NF(G)=[n_{ij}]$ defined by
\[
n_{ij}=\begin{cases}\delta(i)\delta(j), & \text{if } i\sim j \text{ in } G,\\
0, & \text{otherwise.}\end{cases}
\]
\end{definition}

The matrix $\NF(G)$ is real and symmetric; its eigenvalues are denoted by
$\lambda_1\ge\lambda_2\ge\dots\ge\lambda_n$, and $\lambda_1$ is called the \emph{neighbourhood
second Zagreb spectral radius} of $G$, written $\rNF(G)$. Note that the sum of the entries of
$\NF(G)$ above the diagonal is precisely $M_2'(G)$, so that $\NF(G)$ is to the neighbourhood second
Zagreb index what the adjacency matrix is to the number of edges.

Two reasons make the spectral radius of such a matrix worth bounding. The first is internal to
spectral graph theory: for a nonnegative symmetric matrix the spectral radius governs the growth of
the entries of its powers, and hence the behaviour of weighted walk counts; it is also the quantity
that controls the epidemic threshold in the standard susceptible--infected--susceptible model on a
network, where the critical infection rate is the reciprocal of the spectral radius of the
underlying weighted adjacency matrix \cite{VOK09}. Passing from $A(G)$ to $\NF(G)$ amounts to
weighting each contact by the product of the neighbourhood degree sums of its endpoints, so that
$\rNF(G)$ measures propagation in a network in which the transmissibility of a contact grows with
the connectivity of the surroundings of both of its ends. The second reason is that $\rNF(G)$
bounds the neighbourhood second Zagreb energy from both sides, as Theorems~\ref{thm:KM},
\ref{thm:energy-lower} and \ref{thm:energy-T} show, and the energy of the neighbourhood  first Zagreb matrix was shown in
\cite{MBDP22}, by a regression analysis of physico-chemical properties, to be a useful molecular descriptor. Every sharpening of a bound on $\rNF(G)$ therefore tightens the range within which
such a descriptor can lie. We do not undertake a correlation study here; the QSPR analysis of
$\NF$-based descriptors, in the manner of \cite{MBDP22}, would be a natural complement to the
present work.

In \cite{BV25}, sharp bounds for $\rNF(G)$ and for the neighbourhood second Zagreb energy were
obtained, together with the corresponding extremal graphs. Our aim here is twofold. First, we
sharpen several of the estimates available for $\rNF(G)$. Lemma~\ref{lem:Ww} replaces the factor
$\demax^{2}$ in the classical comparison $\rNF(G)\le\demax^{2}\rho(G)$ by the largest edge weight
$W(G)=\max_{\{i,j\}\in E(G)}\delta(i)\delta(j)$. Theorem~\ref{thm:rho5} proves that
$\rNF(G)\ge\rho(G)^{5}$, which is never weaker than the lower estimate $\demin^{2}\rho(G)$, and strictly better unless $G$ is regular or bipartite semiregular. Also this is
the key to the equality case of the Berman-Zhang-type bound of Theorem~\ref{thm:BZ}. Lemma~\ref{lem:delta-const}
allows the equality conditions in Theorem~\ref{thm:delta-rho} and Corollary~\ref{cor:delta52} to be
stated as characterisations rather than as sufficient conditions; the same argument sharpens the
corresponding statements for the matrix $\Upsilon(G)$ of \cite{MBDP22}, see
Remarks~\ref{rem:Ups-delta} and~\ref{rem:Ups-Hong}. Secondly, we determine bounds for classes of
graphs that have not been considered from this point of view, most notably the unicyclic graphs,
for which we show that the cycle $C_{n}$ is the unique minimiser of $\rNF$ on $\mathcal{U}_n$ and that the triangle with
$n-3$ pendant edges attached at one vertex is the unique maximiser.

\subsection*{Contributions and relation to earlier work}

Since several of the estimates below are of a kind that has already been considered for $\NF(G)$
in \cite{BV25}, for the additive matrix $\Upsilon(G)$ in \cite{BV24, MBDP22}, or for the adjacency
matrix itself, we state result by result what the present paper contributes. Table~\ref{tab:new}
summarises this.

\begin{table}[ht]
\centering
\caption{What is new in this paper? The classical results used as tools are
Lemma~\ref{lem:FMS} \cite{FMS93}, Lemma~\ref{lem:Hofmeister} \cite{Hof88}, Hong's inequality
\cite{Hong88}, Nikiforov's inequality \cite{Nik02,Nik09}, Kober's inequality \cite{Kob58},
and the theorems of Smith \cite{Smi70} and of Cheng and Liu \cite{CL07}.}
\label{tab:new}
\small
\begin{tabular}{>{\raggedright\arraybackslash}p{2.6cm}>{\raggedright\arraybackslash}p{5.4cm}>{\raggedright\arraybackslash}p{6.6cm}}
\toprule
Result & Statement & Contribution\\
\midrule
Prop.~\ref{prop:factor} & $\NF=D_\delta AD_\delta$ & elementary; the case of constant $\delta$ occurs in \cite[Prop.~1]{BV25}\\
Prop.~\ref{prop:congruent}, Cor.~\ref{cor:inertia} & $\NF$ congruent to $A$ & new (det formula: \cite[Cor.~1]{BV25}; Cor.~\ref{cor:inertia}(iv): \cite[Prop.~1]{BV25})\\
Lem.~\ref{lem:const-edge} & constant edge weight $\Rightarrow$ $\delta$ constant & new\\
Thm.~\ref{thm:dmax5} & $\rNF\le\dmax^{5}$ & implied by \cite[Cor.~4]{BV25}; dominated\\
Thm.~\ref{thm:rho5} & $\rNF\ge\rho^{5}$ & new\\
Lem.~\ref{lem:Rconst} & constant row sums $\Rightarrow$ $R_0\le\rho^{5}$ & new\\
Thm.~\ref{thm:BZ} & $\rNF\le\max_{i\sim j}\sqrt{R_iR_j}$ & new; equality characterised.\\
Thm.~\ref{thm:delta-rho} & $\demin^{2}\rho\le\rNF\le\demax^{2}\rho$ & \cite[Thm.~9]{BV25}; equality condition corrected.\\
Cor.~\ref{cor:delta52} & $\rNF\le\demax^{5/2}$ & from Thm.~\ref{thm:delta-rho} and \cite{FMS93}; equality characterised.\\
Cor.~\ref{cor:Hong} & (i) $\rNF\le\demax^{2}\sqrt{2m-n+1}$\newline(ii) $\rNF\le\demax^{2}\sqrt{2m(1-1/\omega)}$ & part (i) is \cite[Cor.~3]{BV25}, equality condition corrected; part (ii) and the combined form are new.\\
Lem.~\ref{lem:Ww} & $W\le\rNF$, $w\rho\le\rNF\le W\rho$ & new\\
Thm.~\ref{thm:bipartite}, \ref{thm:trees} & bipartite graphs and trees & new \\
Lem.~\ref{lem:unicyclic-relocation}, \ref{lem:unicyclic-weighted}, Thm.~\ref{thm:unicyclic} & unicyclic graphs: $32=\rNF(C_n)\le\rNF\le\rNF(\Sp n)$, both extremal graphs & new, both equalities are characterised. \\
Thm.~\ref{thm:M1} & $\rNF\ge w\sqrt{M_1/n}$ & new; weaker form \cite[Cor.~4]{BV25}, equality corrected\\
Cor.~\ref{cor:M1-rho5} & $\rNF\ge(M_1/n)^{5/2}$ & new\\
Thm.~\ref{thm:KM} & energy, upper bounds &  new\\
Thm.~\ref{thm:energy-lower}, \ref{thm:energy-T}, Cor.~\ref{cor:unicyclic-energy} & energy, lower bounds & new for $\NF$ (standard arguments); extremal graphs characterised\\
\bottomrule
\end{tabular}
\end{table}


The paper is organised as follows. Section~\ref{sec:prelim} collects the matrix tools that will be
used and records the value of $\rNF$ for several basic families of graphs.
Section~\ref{sec:general} contains bounds valid for all connected graphs; in particular
Theorem~\ref{thm:rho5} gives the lower bound $\rNF(G)\ge\rho(G)^{5}$, Theorem~\ref{thm:BZ} gives a
Berman--Zhang-type upper bound whose extremal graphs are exactly the regular and bipartite
semiregular graphs, and Lemma~\ref{lem:Ww} improves the classical estimate
$\rNF(G)\le\demax^{2}\rho(G)$ by replacing $\demax^{2}$ with the largest edge weight of $\NF(G)$.
Section~\ref{sec:bipartite} treats bipartite graphs and trees, and Section~\ref{sec:unicyclic} is
devoted to unicyclic graphs, where a two-sided bound with both extremal cases discussed is
obtained. Section~\ref{sec:M1} relates $\rNF(G)$ to the first Zagreb index and
Section~\ref{sec:energy} contains bounds for the neighbourhood second Zagreb energy.

\section{Preliminaries}\label{sec:prelim}

We use the following standard facts; see \cite[Chapters~1, 4 and 8]{HJ13} and \cite{CRS10}.

\begin{lemma}[Rayleigh--Ritz {\cite{HJ13}}]\label{lem:RR} 
Let $M$ be a real symmetric $n\times n$ matrix with the largest eigenvalue $\lambda_{\max}(M)$. Then
$\lambda_{\max}(M)=\max_{\|x\|=1}x^{T}Mx$, and $x^{T}Mx\le\lambda_{\max}(M)$ for every unit vector
$x$, with equality if and only if $x$ is an eigenvector of $M$ associated with $\lambda_{\max}(M)$.
\end{lemma}

\begin{lemma}[Perron--Frobenius {\cite{HJ13}}]\label{lem:PF} 
Let $B=[b_{ij}]$ and $C=[c_{ij}]$ be nonnegative $n\times n$ matrices.
\begin{enumerate}[label=(\alph*)]
\item If $b_{ij}\le c_{ij}$ for all $i,j$ and $C$ is irreducible, then $\rho(B)\le\rho(C)$, with
equality if and only if $B=C$.
\item If $C$ is irreducible with row sums $R_1,\dots,R_n$, then
$\min_iR_i\le\rho(C)\le\max_iR_i$, and either equality holds if and only if $R_1=\dots=R_n$.
\item \emph{(Collatz--Wielandt)} If $C$ is irreducible and $z\in\mathbb{R}^{n}$ is a positive vector,
then $\rho(C)\le\max_i(Cz)_i/z_i$.
\end{enumerate}
\end{lemma}

\begin{lemma}[Cauchy interlacing {\cite{HJ13}}]\label{lem:interlacing} 
If $M$ is a real symmetric $n\times n$ matrix and $M'$ is a principal submatrix of $M$, then
$\lambda_{\max}(M')\le\lambda_{\max}(M)$.
\end{lemma}

If $G$ is connected and $n\ge 2$, then $A(G)$ is irreducible and all the weights $\delta(i)\delta(j)$,
$\{i,j\}\in E(G)$, are positive; hence $\NF(G)$ is a nonnegative irreducible matrix and $\rNF(G)$
is a simple eigenvalue admitting a positive eigenvector, the \emph{Perron vector} of $\NF(G)$. The
next observation, although elementary, explains the structure of $\NF(G)$ and will be used
repeatedly.

\begin{proposition}\label{prop:factor}
Let $G$ be a graph with no isolated vertices and let $D_\delta=\diag\bigl(\delta(1),\dots,\delta(n)\bigr)$.
Then
\[
\NF(G)=D_\delta\,A(G)\,D_\delta,
\]
and consequently $\rNF(G)=\rho\bigl(A(G)D_\delta^{2}\bigr)$. In particular, if $\delta(i)=c$ for
every $i\in V(G)$, then $\NF(G)=c^{2}A(G)$ and $\rNF(G)=c^{2}\rho(G)$.
\end{proposition}

\begin{proof}
The $(i,j)$-th entry of $D_\delta A(G)D_\delta$ equals $\delta(i)a_{ij}\delta(j)$, which is
$\delta(i)\delta(j)$ if $i\sim j$ and $0$ otherwise. The matrices $D_\delta A(G)D_\delta$ and
$A(G)D_\delta^{2}$ are similar via $D_\delta$, which is invertible because $\delta(i)>0$ for every
$i$ when $G$ has no isolated vertex.
\end{proof}

The factorisation in Proposition~\ref{prop:factor} is a congruence, and this has consequences that
are not available for the additive matrix $\Upsilon(G)=[\delta(i)+\delta(j)]$ of \cite{MBDP22}.

\begin{proposition}\label{prop:congruent}
Let $G$ be a graph with no isolated vertices. Then $\NF(G)$ and $A(G)$ are congruent; consequently
they have the same inertia, that is, the same number of positive, negative and zero eigenvalues.
In particular
\[
\rank\NF(G)=\rank A(G)\qquad\text{and}\qquad\det\NF(G)=\Bigl(\prod_{i=1}^{n}\delta(i)\Bigr)^{2}\det A(G).
\]
\end{proposition}

\begin{proof}
Since $G$ has no isolated vertex, $\delta(i)>0$ for every $i$, so $D_\delta$ is nonsingular and
$\NF(G)=D_\delta A(G)D_\delta$ exhibits a congruence between $\NF(G)$ and $A(G)$. Sylvester's law
of inertia \cite{HJ13} gives the equality of inertias, and hence of ranks. The determinant formula
follows from $\det\NF(G)=(\det D_\delta)^{2}\det A(G)$.
\end{proof}

Proposition~\ref{prop:congruent} is a transfer principle: every statement about the inertia of
$A(G)$ holds verbatim for $\NF(G)$. We record the consequences used later. Note that each of them
would require a separate argument for the additive matrix $\Upsilon(G)$, for which no analogue of
Proposition~\ref{prop:factor} is available.
\begin{corollary}[part (iv): {\cite [Prop.~1(i)]{BV25}}]\label{cor:inertia}
Let $G$ be a graph on $n$ vertices with no isolated vertices. Then
\begin{enumerate}[label=(\roman*)]
\item the nullity of $\NF(G)$ equals the nullity of $G$; in particular $\NF(G)$ is nonsingular if
and only if $A(G)$ is;
\item $\rank\NF(G)=2$ if and only if $G\cong K_{p,q}$ for some $p,q\ge1$;
\item $\NF(G)$ has exactly one positive eigenvalue if and only if $G$ is a complete multipartite
graph;
\item if $G$ is bipartite, then the spectrum of $\NF(G)$ is symmetric about the origin.
\end{enumerate}
\end{corollary}

\begin{proof}
Parts (i)--(iii) follow from Proposition~\ref{prop:congruent} together with, respectively, the
definition of nullity, the theorem of Cheng and Liu \cite{CL07} describing the graphs of nullity
$n-2$ (these are the complete bipartite graphs together with isolated vertices, and there are no
isolated vertices here), and the theorem of Smith \cite{Smi70} characterising the graphs with
exactly one positive adjacency eigenvalue. For (iv), let $(X,Y)$ be a bipartition of $G$ and let
$J=\diag(\varepsilon_1,\dots,\varepsilon_n)$ with $\varepsilon_i=1$ for $i\in X$ and
$\varepsilon_i=-1$ for $i\in Y$. Every edge of $G$ joins $X$ to $Y$, so $J\NF(G)J=-\NF(G)$, and
$J$ is orthogonal; hence the spectrum of $\NF(G)$ is symmetric.
\end{proof}

For a graph $G$ put
\begin{equation}\label{eq:T}
T(G)=\sum_{\{i,j\}\in E(G)}\delta(i)^{2}\delta(j)^{2},
\end{equation}

\begin{lemma}\label{lem:moments}
Let $G$ be a graph on $n$ vertices and let $\lambda_1\ge\dots\ge\lambda_n$ be the eigenvalues of
$\NF(G)$. Then
\[
\sum_{i=1}^{n}\lambda_i=0,\qquad\sum_{i=1}^{n}\lambda_i^{2}=\tr\bigl(\NF(G)^{2}\bigr)=2T(G),\qquad
\sum_{i<j}\lambda_i\lambda_j=-T(G).
\]
\end{lemma}

\begin{proof}
The diagonal of $\NF(G)$ is zero, whence $\sum_i\lambda_i=\tr\NF(G)=0$. Since $\NF(G)$ is
symmetric, $\tr(\NF(G)^{2})=\sum_{i,j}n_{ij}^{2}=2\sum_{\{i,j\}\in E(G)}\delta(i)^{2}\delta(j)^{2}=2T(G)$.
Finally $2\sum_{i<j}\lambda_i\lambda_j=\bigl(\sum_i\lambda_i\bigr)^{2}-\sum_i\lambda_i^{2}=-2T(G)$.
\end{proof}

\begin{proposition}\label{prop:families}
The following hold.
\begin{enumerate}[label=(\roman*)]
\item $\rNF(K_n)=(n-1)^{5}$ for $n\ge2$.
\item $\rNF(K_{p,q})=(pq)^{2}\sqrt{pq}$ for $p,q\ge1$; in particular
$\rNF(K_{1,n-1})=(n-1)^{2}\sqrt{n-1}$.
\item $\rNF(C_n)=32$ for every $n\ge3$.
\item If $G$ is $r$-regular, then $\rNF(G)=r^{4}\rho(G)=r^{5}$.
\end{enumerate}
\end{proposition}

\begin{proof}
In each case $\delta$ is constant, so Proposition~\ref{prop:factor} applies. For $K_n$ we have
$\delta\equiv(n-1)^{2}$ and $\rho(K_n)=n-1$; for $K_{p,q}$ every vertex satisfies $\delta\equiv pq$
and $\rho(K_{p,q})=\sqrt{pq}$; for $C_n$ we have $\delta\equiv4$ and $\rho(C_n)=2$, whence
$\rNF(C_n)=16\cdot2=32$; for an $r$-regular graph $\delta\equiv r^{2}$ and $\rho(G)=r$.
\end{proof}

Note that (iii) shows that $\rNF(C_n)$ does not depend on $n$: all cycles have the same
neighbourhood second Zagreb spectral radius.

\section{Bounds for general graphs}\label{sec:general}

We begin with three lemmas describing when $\delta$ is constant. Recall that a bipartite graph
with bipartition $(X,Y)$ is \emph{$(r,s)$-semiregular} if every vertex of $X$ has degree $r$ and
every vertex of $Y$ has degree $s$.

\begin{lemma}\label{lem:delta-const}
Let $G$ be a connected graph on $n$ vertices. Then $\delta$ is constant on $ V(G)$  $(\mbox{i.e.}~ \demin=\demax)$ 
if and only if $G$ is regular or bipartite semiregular.
\end{lemma}

\begin{proof}
Assume $\delta(i)=p$ for all $i\in V(G)$. Choose $k,t\in V(G)$ with $d(k)=\dmin$ and $d(t)=\dmax$.
Since every neighbour of $k$ has degree at most $\dmax$,
\[
p=\delta(k)=\sum_{j\sim k}d(j)\le d(k)\dmax=\dmin\dmax,
\]
and since every neighbour of $t$ has degree at least $\dmin$,
\[
p=\delta(t)=\sum_{j\sim t}d(j)\ge d(t)\dmin=\dmax\dmin.
\]
Hence $p=\dmin\dmax$. If $\dmin=\dmax$, then $G$ is regular. Assume therefore $\dmin\ne\dmax$.
From $\delta(i)=\dmin\dmax$ for all $i$ we deduce: if $d(i)=\dmin$, then all neighbours of $i$
have degree $\dmax$; and if $d(i)=\dmax$, then all neighbours of $i$ have degree $\dmin$. Indeed,
in the first case $\dmin\dmax=\delta(i)=\sum_{j\sim i}d(j)\le\dmin\dmax$ forces $d(j)=\dmax$ for
every $j\sim i$, and the second case is symmetric.

Fix $k$ with $d(k)=\dmin$ and put $V_1=\{k\}$, $V_l=N_G(V_{l-1})$ for $l\ge2$ and $V(G)=\bigcup_{l\ge1} V_{l}$. By the above, the
vertices of $V_1,V_3,V_5,\dots$ have degree $\dmin$ and those of $V_2,V_4,\dots$ have degree
$\dmax$, and every edge joins a vertex of degree $\dmin$ to a vertex of degree $\dmax$. Since $G$
is connected, $V(G)=\bigcup_{m\ge1}V_m$, and therefore $V(G)$ is partitioned into the set of
vertices of degree $\dmin$ and the set of vertices of degree $\dmax$, all edges joining the two
sets. Thus $G$ is $(\dmin,\dmax)$-semiregular bipartite.

Conversely, if $G$ is $r$-regular, then $\delta(i)=r^{2}$ for every $i$; and if $G$ is
$(r,s)$-semiregular bipartite, then every vertex is adjacent only to vertices of the opposite part,
so $\delta(i)=rs$ for every $i$.
\end{proof}

\begin{lemma}[Favaron et al. \cite{FMS93}, Cao \cite{Cao98}]\label{lem:FMS}
Let $G$ be a connected graph on $n$ vertices. Then $\rho(G)\le\sqrt{\demax}$, with equality if and
only if $G$ is regular or bipartite semiregular.
\end{lemma}

We include the short proof for completeness.

\begin{proof}
Let $X=(x_{1}, \cdots, x_{n})^{T}$ be the positive unit eigenvector (Perron vector) of $A(G)$ corresponding to $\rho(G)$. For every
$i$, the Cauchy--Schwarz inequality gives
\[
\rho(G)^{2}x_i^{2}=\Bigl(\sum_{j\sim i}x_j\Bigr)^{2}\le d(i)\sum_{j\sim i}x_j^{2}.
\]
Summing over $i$ and interchanging the order of summation,
\[
\rho(G)^{2}=\rho(G)^{2}\sum_{i=1}^{n}x_i^{2}\le\sum_{i=1}^{n}d(i)\sum_{j\sim i}x_j^{2}
=\sum_{j=1}^{n}x_j^{2}\sum_{i\sim j}d(i)=\sum_{j=1}^{n}x_j^{2}\,\delta(j)\le\demax.
\]
Equality forces $\delta(j)=\demax$ for every $j$ (the Perron vector is positive), and
Lemma~\ref{lem:delta-const} completes the proof. Conversely, if $\delta$ is constant, say
$\delta\equiv c$, then $A^{2}\one=c\one$; since $\one > 0, c= \rho(A^2)=\rho(G)^{2}$ by the Perron-Frobenius theorem.
\end{proof}

Constancy of the edge weights already forces constancy of $\delta$; this is what turns several
equality conditions below into structural characterisations.

\begin{lemma}\label{lem:const-edge}
Let $G$ be a connected graph and suppose that $\delta(i)\delta(j)$ takes one and the same value on
every edge of $G$. Then $\delta$ is constant on $V(G)$, and hence, by Lemma~\ref{lem:delta-const},
$G$ is regular or bipartite semiregular.
\end{lemma}

\begin{proof}
Let $\delta(i)\delta(j)=c>0$ for every $\{i,j\}\in E(G)$. If two vertices $i$ and $k$ have a common
neighbour $j$, then $\delta(i)\delta(j)=c=\delta(k)\delta(j)$ and therefore $\delta(i)=\delta(k)$.
Thus $\delta$ is constant on each class of the relation ``joined by a walk of even length''. If $G$
is connected and not bipartite, this relation has a single class and $\delta$ is constant.

Suppose then that $G$ is connected and bipartite with bipartition $(X,Y)$, so that $\delta\equiv a$
on $X$ and $\delta\equiv b$ on $Y$ for some $a,b>0$. Every edge of $G$ has exactly one endpoint in
$X$, whence
\[
\sum_{i\in X}d(i)\delta(i)=\sum_{i\in X}\sum_{j\sim i}d(i)d(j)=\sum_{\{i,j\}\in E(G)}d(i)d(j)=M_2(G),
\]
$M_2$ denoting the second Zagreb index, and the same computation over $Y$ gives
$\sum_{j\in Y}d(j)\delta(j)=M_2(G)$. Consequently $a\sum_{i\in X}d(i)=b\sum_{j\in Y}d(j)$, and since
$\sum_{i\in X}d(i)=\sum_{j\in Y}d(j)=m>0$ we conclude $a=b$.
\end{proof}

\begin{theorem}[cf.\ {\cite[Cor.~4]{BV25}}]\label{thm:dmax5}
Let $G$ be a connected graph on $n$ vertices. Then
\[
\rNF(G)\le\dmax^{5},
\]
and equality holds if and only if $G$ is regular.
\end{theorem}

\begin{proof}
By \eqref{eq:delta-range}, $\delta(i)\le\dmax^{2}$ for every $i$. The $i$-th row sum of $\NF(G)$
equals
\begin{equation}\label{eq:rowsum}
R_i=\delta(i)\sum_{j\sim i}\delta(j)\le\dmax^{2}\cdot d(i)\dmax^{2}\le\dmax^{5}.
\end{equation}
Since $\NF(G)$ is nonnegative and irreducible, Lemma~\ref{lem:PF} yields
$\rNF(G)\le\max_iR_i\le\dmax^{5}$.

If $\rNF(G)=\dmax^{5}$, then $\max_iR_i=\dmax^{5}=\rNF(G)$, so by Lemma~\ref{lem:PF} all row sums
are equal to $\dmax^{5}$. Hence, both inequalities in \eqref{eq:rowsum} are equalities for every
$i$, which forces $\delta(i)=\dmax^{2}$ and $d(i)=\dmax$ for every $i\in V(G)$; that is, $G$ is
regular. Conversely, if $G$ is $\dmax$-regular, then $\rNF(G)=\dmax^{5}$ by
Proposition~\ref{prop:families}(iv).
\end{proof}

\begin{remark}\label{rem:n-only}
Since $\dmax\le n-1$, 
Theorem~\ref{thm:dmax5} also gives $\rNF(G)\le(n-1)^{5}$, with equality if and only if $G \cong K_n$. For sparse
graphs the degree-based bound is far better: for a graph on $n=100$ vertices with $\dmax=10$ one
gets $\rNF(G)\le10^{5}$, while the bound in terms of $n$ only gives $99^{5}$.
\end{remark}

Theorem~\ref{thm:dmax5} is an upper estimate. We now give a lower estimate in terms of the adjacency
spectral radius alone; besides its own interest, it is the key to the equality case of
Theorem~\ref{thm:BZ} below.

\begin{theorem}\label{thm:rho5}
Let $G$ be a connected graph on $n\ge2$ vertices. Then
\[
\rNF(G)\ge\rho(G)^{5},
\]
and equality holds if and only if $G$ is regular or bipartite semiregular.
\end{theorem}

\begin{proof}
Write $\rho=\rho(G)>0$ and $A=A(G)$, let $\bxi=(\xi_1,\dots,\xi_n)^{T}$ be the positive unit
eigenvector of $A$ corresponding to $\rho$, and put $Q=A^{2}$. The entry $Q_{ij}$ is the number of
walks of length $2$ from $i$ to $j$, so $Q$ is nonnegative and symmetric and
$\sum_{j}Q_{ij}=\sum_{k\sim i}d(k)=\delta(i)$. Moreover $\delta(i)>0$ for every $i$, so
$D_\delta=\diag(\delta(1),\dots,\delta(n))$ is invertible.

\emph{Step 1: $\sum_{i}\xi_i^{2}/\delta(i)\le1/\rho^{2}$.} For every $x\in\mathbb{R}^{n}$,
\[
x^{T}(D_\delta-Q)x=\sum_{i}\delta(i)x_i^{2}-\sum_{i,j}Q_{ij}x_ix_j
=\tfrac12\sum_{i,j}Q_{ij}(x_i-x_j)^{2}\ge0,
\]
and $x^{T}Qx=\|Ax\|^{2}\ge0$. Hence the symmetric matrix
$\Gamma=D_\delta^{-1/2}QD_\delta^{-1/2}$ satisfies $0\preceq\Gamma\preceq I$, so all its eigenvalues
lie in $[0,1]$ and $\Gamma^{2}\preceq\Gamma$. Put $q=D_\delta^{1/2}\bxi$. Then
$\Gamma q=D_\delta^{-1/2}Q\bxi=\rho^{2}D_\delta^{-1/2}\bxi$ and $q^{T}\Gamma q=\bxi^{T}Q\bxi=\rho^{2}$,
so that
\[
\rho^{4}\sum_{i}\frac{\xi_i^{2}}{\delta(i)}=\|\Gamma q\|^{2}=q^{T}\Gamma^{2}q\le q^{T}\Gamma q=\rho^{2}.
\]

\emph{Step 2.} Put $\ell=\sum_{i}\xi_i^{2}\log\delta(i)$. Since $\sum_i\xi_i^{2}=1$, Jensen's
inequality for the concave function $\log$ and Step~1 give
\[
-\ell=\sum_{i}\xi_i^{2}\log\frac{1}{\delta(i)}\le\log\sum_{i}\frac{\xi_i^{2}}{\delta(i)}\le-\log\rho^{2},
\qquad\text{that is,}\qquad e^{\ell}\ge\rho^{2}.
\]

\emph{Step 3: $\bxi^{T}\NF\bxi\ge\rho\,e^{2\ell}$.} For every edge $\{i,j\}$ put
$x_{ij}=\log\delta(i)+\log\delta(j)-2\ell$, so that
$\delta(i)\delta(j)e^{-2\ell}=e^{x_{ij}}\ge1+x_{ij}$, with equality if and only if $x_{ij}=0$.
Multiplying by $a_{ij}\xi_i\xi_j$, which is positive when $i\sim j$, and summing over all ordered
pairs $(i,j)$, we obtain
\begin{align*}
e^{-2\ell}\,\bxi^{T}\NF\bxi
&\ge\sum_{i,j}a_{ij}\xi_i\xi_j
+\sum_{i,j}a_{ij}\xi_i\xi_j\bigl[(\log\delta(i)-\ell)+(\log\delta(j)-\ell)\bigr]\\
&=\rho+2\rho\sum_{i}\xi_i^{2}\bigl(\log\delta(i)-\ell\bigr)=\rho,
\end{align*}
because $\sum_{j}a_{ij}\xi_j=\rho\xi_i$ and $\sum_i\xi_i^{2}(\log\delta(i)-\ell)=\ell-\ell=0$.

\emph{Step 4.} By Lemma~\ref{lem:RR}, Step~3 and Step~2,
\[
\rNF(G)\ge\bxi^{T}\NF\bxi\ge\rho\,e^{2\ell}\ge\rho\cdot\rho^{4}=\rho^{5}.
\]

\emph{Equality.} If $\rNF(G)=\rho^{5}$, then equality holds in Step~3, so $x_{ij}=0$ for every
edge, i.e.\ $\delta(i)\delta(j)=e^{2\ell}$ for every $\{i,j\}\in E(G)$. By
Lemma~\ref{lem:const-edge}, $G$ is regular or bipartite semiregular. Conversely, if $G$ is regular
or bipartite semiregular, then $\delta\equiv c$ by Lemma~\ref{lem:delta-const}, $\rho^{2}=c$ by
Lemma~\ref{lem:FMS}, and $\rNF(G)=c^{2}\rho=\rho^{5}$ by Proposition~\ref{prop:factor}.
\end{proof}

Theorem~\ref{thm:dmax5} is a row-sum estimate. Weighting the rows against each other gives a bound
of the Berman--Zhang type \cite{BZ01} which is never weaker (Remark~\ref{rem:BZ}). For its equality
case we need the following estimate for graphs whose row sums $R_i$ are constant on the classes of
vertices joined by walks of even length.

\begin{lemma}\label{lem:Rconst}
Let $G$ be a connected graph on $n\ge2$ vertices and let $R_i=\delta(i)\sum_{j\sim i}\delta(j)$ be
the $i$-th row sum of $\NF(G)$.
\begin{enumerate}[label=(\roman*)]
\item If $R_i=R_0$ for every $i\in V(G)$, then $R_0\le\rho(G)^{5}$.
\item If $G$ is bipartite with bipartition $(X,Y)$, and $R_i=R_X$ for $i\in X$ and $R_i=R_Y$ for
$i\in Y$, then $R_XR_Y\le\rho(G)^{10}$.
\end{enumerate}
\end{lemma}

\begin{proof}
Write $\rho=\rho(G)$, $A=A(G)$ and $\sigma(i)=\sum_{j\sim i}\delta(j)$, so that
$R_i=\delta(i)\sigma(i)$; in vector form $\boldsymbol{\delta}=(\delta(1), \cdots ,\delta(n))^{T}=A^{2}\one$ and $\boldsymbol{\sigma}=A \boldsymbol{\delta}$. Let $\bxi>0$ be an
eigenvector of $A$ belonging to $\rho$.

(i) We have $\bxi^{T}\boldsymbol{\delta}=(A^{2}\bxi)^{T}\one=\rho^{2}\sum_i\xi_i$ and
$\bxi^{T}\boldsymbol{\sigma}=(A\bxi)^{T}\boldsymbol{\delta}=\rho\,\bxi^{T}\boldsymbol{\delta}=\rho^{3}\sum_i\xi_i$. Since
$\sigma(i)=R_0/\delta(i)$, this gives $\sum_i\xi_i/\delta(i)=(\rho^{3}/R_0)\sum_i\xi_i$. By the
Cauchy--Schwarz inequality,
\[
\Bigl(\sum_i\xi_i\Bigr)^{2}=\Bigl(\sum_i\sqrt{\xi_i\delta(i)}\sqrt{\xi_i/\delta(i)}\Bigr)^{2}
\le\Bigl(\sum_i\xi_i\delta(i)\Bigr)\Bigl(\sum_i\frac{\xi_i}{\delta(i)}\Bigr)
=\frac{\rho^{5}}{R_0}\Bigl(\sum_i\xi_i\Bigr)^{2}.
\]
So $R_0\le\rho^{5}$.

(ii) Let $\bxi_X$ and $\bxi_Y$ be the vectors that agree with $\bxi$ on $X$ and on $Y$,
respectively, and vanish elsewhere. Since every edge joins $X$ to $Y$, we have
$A\bxi_X=\rho\bxi_Y$ and $A\bxi_Y=\rho\bxi_X$, hence $A^{2}\bxi_X=\rho^{2}\bxi_X$ and
$A^{2}\bxi_Y=\rho^{2}\bxi_Y$. Therefore
\[
\sum_{i\in X}\xi_i\delta(i)=(A^{2}\bxi_X)^{T}\one=\rho^{2}\sum_{i\in X}\xi_i,\qquad
\sum_{i\in X}\xi_i\sigma(i)=(A\bxi_X)^{T}\boldsymbol{\delta}=\rho\sum_{j\in Y}\xi_j\delta(j)=\rho^{3}\sum_{j\in Y}\xi_j.
\]
As $\sigma(i)=R_X/\delta(i)$ for $i\in X$, the Cauchy--Schwarz inequality over $X$ gives
\[
\Bigl(\sum_{i\in X}\xi_i\Bigr)^{2}\le\Bigl(\sum_{i\in X}\xi_i\delta(i)\Bigr)\Bigl(\sum_{i\in X}\frac{\xi_i}{\delta(i)}\Bigr)
=\rho^{2}\sum_{i\in X}\xi_i\cdot\frac{\rho^{3}}{R_X}\sum_{j\in Y}\xi_j,
\]
that is, $R_X\sum_{i\in X}\xi_i\le\rho^{5}\sum_{j\in Y}\xi_j$. By symmetry,
$R_Y\sum_{j\in Y}\xi_j\le\rho^{5}\sum_{i\in X}\xi_i$. Multiplying these two inequalities between positive numbers gives $R_XR_Y\le\rho^{10}$.
\end{proof}

\begin{theorem}\label{thm:BZ}
Let $G$ be a connected graph on $n\ge2$ vertices and, for $i\in V(G)$, let
$R_i=\delta(i)\sum_{j\sim i}\delta(j)$ be the $i$-th row sum of $\NF(G)$. Put
\[
\beta(G)=\max_{\{i,j\}\in E(G)}\sqrt{R_iR_j}.
\]
Then $\rNF(G)\le\beta(G)$. Moreover, equality holds if and only if $G$ is regular or $G$ is bipartite semiregular.
\end{theorem}

\begin{proof}
Let $z=(\sqrt{R_1},\dots,\sqrt{R_n})^{T}>0$ and let $\bpi>0$ be a Perron vector of $\NF(G)$, so
that $\NF\bpi=\rNF(G)\bpi$.

\emph{Step 1: the inequality.} For every $i$,
\begin{equation}\label{eq:BZ-chain}
(\NF z)_i=\sum_{j\sim i}\delta(i)\delta(j)\sqrt{R_j}
\;\overset{\mathrm{(I)}}{\le}\;R_i\max_{j\sim i}\sqrt{R_j}
\;\overset{\mathrm{(II)}}{\le}\;\beta(G)\sqrt{R_i},
\end{equation}
where (I) uses $\sum_{j\sim i}\delta(i)\delta(j)=R_i$, and (II) holds because
$\sqrt{R_i}\max_{j\sim i}\sqrt{R_j}$ is one of the numbers $\sqrt{R_iR_j}$, $j\sim i$. Thus
$\NF z\le\beta(G)z$ entrywise. Since $\NF$ is symmetric,
\begin{equation}\label{eq:BZ-perron}
\rNF(G)\,\bpi^{T}z=\bpi^{T}\NF z\le\beta(G)\,\bpi^{T}z,
\end{equation}
and $\bpi^{T}z>0$ gives $\rNF(G)\le\beta(G)$.

\emph{Step 2: regular and bipartite semiregular graphs attain equality.} If $G$ is regular or
bipartite semiregular, then $\delta\equiv c$ by Lemma~\ref{lem:delta-const}, and $d(i)d(j)=c$ on
every edge ($r\cdot r=r^{2}$ for an $r$-regular graph, $rs$ for an $(r,s)$-semiregular bipartite
graph). Hence $R_i=c\cdot c\,d(i)=c^{2}d(i)$ and $\sqrt{R_iR_j}=c^{2}\sqrt{d(i)d(j)}=c^{5/2}$ for
every edge, so $\beta(G)=c^{5/2}$. By Proposition~\ref{prop:factor} and Lemma~\ref{lem:FMS},
$\rNF(G)=c^{2}\rho(G)=c^{2}\sqrt{c}=\beta(G)$.

\emph{Step 3: equality forces $R$ to be constant on the classes of even walks.} Suppose
$\rNF(G)=\beta(G)$. Then equality holds in \eqref{eq:BZ-perron}, i.e.\
$\bpi^{T}\bigl(\beta(G)z-\NF z\bigr)=0$; as $\bpi>0$ and $\beta(G)z-\NF z\ge0$, we get
$\NF z=\beta(G)z$, so (I) in \eqref{eq:BZ-chain} is an equality for every $i$. The weights
$\delta(i)\delta(j)$ are positive, and a weighted average equals the maximum of the averaged values
only if all of them are equal; hence all neighbours of each vertex have the same value of $R$. If
$j$ is a common neighbour of $i$ and $k$, then $i$ and $k$ are both neighbours of $j$, so $R_i=R_k$.
Consequently $R$ is constant on every class of vertices joined by walks of even length. If $G$ is
not bipartite, any two vertices are joined by an even walk, so $R\equiv R_0$ on $V(G)$ and
$\beta(G)=R_0$. If $G$ is bipartite with bipartition $(X,Y)$, then $R\equiv R_X$ on $X$ and
$R\equiv R_Y$ on $Y$, every edge joins $X$ to $Y$, and $\beta(G)=\sqrt{R_XR_Y}$.

\emph{Step 4: equality forces $G$ to be regular or bipartite semiregular.} In either case of
Step~3, Lemma~\ref{lem:Rconst} gives $\beta(G)\le\rho(G)^{5}$. Together with
Theorem~\ref{thm:rho5},
\[
\rho(G)^{5}\le\rNF(G)=\beta(G)\le\rho(G)^{5},
\]
so $\rNF(G)=\rho(G)^{5}$, and the equality part of Theorem~\ref{thm:rho5} shows that $G$ is regular
or bipartite semiregular. Together with Step~2 this proves the equivalence.
\end{proof}

\begin{remark}\label{rem:BZ}
(a) Since $\sqrt{R_iR_j}\le\max_kR_k$ and, by \eqref{eq:rowsum}, $R_k\le\dmax^{5}$,
Theorem~\ref{thm:BZ} yields
\[
\rNF(G)\le\beta(G)\le\max_{i\in V(G)}R_i\le\dmax^{5}.
\]
Thus it refines both the row-sum bound of Lemma~\ref{lem:PF}(b) and Theorem~\ref{thm:dmax5}. In
the regular case all four quantities equal $r^{5}$; for the star $K_{1,n-1}$ one has
$\beta=(n-1)^{5/2}=\rNF$, whereas $\max_iR_i=(n-1)^{3}$ and $\dmax^{5}=(n-1)^{5}$.

(b) Unlike Theorem~\ref{thm:dmax5}, the bound of Theorem~\ref{thm:BZ} is not dominated by
Corollary~\ref{cor:delta52}: for the star $K_{1,n-1}$ both give the exact value $(n-1)^{5/2}$,
whereas for many graphs each is the smaller one. Among the $141$ connected graphs on $3$ to $6$
vertices, $\beta(G)<\demax^{5/2}$  for $100$, $\beta(G)=\demax^{5/2}$ for $17$ and $\beta(G)>\demax^{5/2}$ for  $24$.
\end{remark}

The following theorem relates $\rNF(G)$ and $\rho(G)$. In \cite[Theorem~9]{BV25} equality on either side is asserted
to hold if and only if $G$ is regular. This is incomplete: every bipartite semiregular graph attains equality as well; for instance, for the stars $K_{1,n-1}$ one has $\delta\equiv n-1$ and
$\NF=(n-1)^{2}A$. The corrected statement is as follows.
\begin{theorem}\label{thm:delta-rho}
Let $G$ be a connected graph. Then
\[
\demin^{2}\rho(G)\le\rNF(G)\le\demax^{2}\rho(G).
\]
Equality holds in either inequality if and only if $G$ is regular or bipartite semiregular.
\end{theorem}

\begin{proof}
Let $X=(x_1,\dots,x_n)^{T}$ be the (positive) unit Perron vector of $\NF(G)$. By
Lemma~\ref{lem:RR},
\begin{align*}
\rNF(G)=X^{T}\NF X&=2\sum_{\{i,j\}\in E(G)}\delta(i)\delta(j)x_ix_j\le2\demax^{2}\sum_{\{i,j\}\in E(G)}x_ix_j\\
&=\demax^{2}X^{T}A(G)X\le\demax^{2}\rho(G).
\end{align*}
Since $X>0$, equality in the middle step holds if and only if $\delta(i)\delta(j)=\demax^{2}$ for
every edge $\{i,j\}$, i.e.\ if and only if $\delta\equiv\demax$; by Lemma~\ref{lem:delta-const}
this happens exactly when $G$ is regular or bipartite semiregular.

Similarly, if $Y=(y_1,\dots,y_n)^{T}$ denotes the positive unit Perron vector of $A(G)$, then
\begin{align*}
\rNF(G)\ge Y^{T}\NF Y&=2\sum_{\{i,j\}\in E(G)}\delta(i)\delta(j)y_iy_j\ge2\demin^{2}\sum_{\{i,j\}\in E(G)}y_iy_j\\
&=\demin^{2}Y^{T}A(G)Y=\demin^{2}\rho(G),
\end{align*}
and equality holds if and only if $\delta\equiv\demin$, that is, again if and only if $G$ is
regular or bipartite semiregular.
\end{proof}

\begin{remark}\label{rem:Ups-delta}
Write $\rUps(G)$ for
the spectral radius of the neighbourhood first Zagreb matrix $\Upsilon(G)=[\delta(i)+\delta(j)]$ of
\cite{MBDP22}. The estimate corresponding to Theorem~\ref{thm:delta-rho}, namely
$2\demin\rho(G)\le\rUps(G)\le2\demax\rho(G)$, is \cite[Theorem~6]{MBDP22}, where equality is stated
as holding if $G$ is regular. Since $2\demin A(G)\le\Upsilon(G)\le2\demax A(G)$ entrywise, with
equality on either side exactly when $\delta$ is constant, Lemma~\ref{lem:PF}(a) and
Lemma~\ref{lem:delta-const} show that equality in either inequality holds if and only if $G$ is
regular or bipartite semiregular, so that the bipartite semiregular graphs are extremal there as
well.
\end{remark}

Combining Lemma~\ref{lem:FMS} and Theorem~\ref{thm:delta-rho} we obtain the following bound, which
depends on $\demax$ alone.

\begin{corollary}\label{cor:delta52}
Let $G$ be a connected graph on $n$ vertices. Then
\[
\rNF(G)\le\demax^{2}\sqrt{\demax}=\demax^{5/2},
\]
and equality holds if and only if $G$ is regular or bipartite semiregular.
\end{corollary}

\begin{remark}\label{rem:dominate}
Since $\demax\le\dmax^{2}$, Corollary~\ref{cor:delta52} implies
$\rNF(G)\le\demax^{5/2}\le\dmax^{5}$; thus Corollary~\ref{cor:delta52} is never weaker than
Theorem~\ref{thm:dmax5}, and it is strictly stronger whenever $\demax<\dmax^{2}$, i.e.\ whenever
every vertex of maximum degree has a neighbour of degree less than $\dmax$. (Indeed
$\delta(i)=\dmax^{2}$ holds exactly when $d(i)=\dmax$ and all neighbours of $i$ have degree
$\dmax$.)
\end{remark}

Combining Theorem~\ref{thm:delta-rho} with Hong's inequality \cite{Hong88} and with Nikiforov's
inequality \cite{Nik02, Nik09} gives bounds in terms of the order, the size, and the clique number.
A \emph{complete regular $\omega$-partite graph} is a complete $\omega$-partite graph
$K_{t,t,\dots,t}$ whose $\omega$ parts all have the same size $t$; for $t=1$ it is $K_\omega$.

\begin{corollary}[part (i): {\cite[Cor.~3]{BV25}}]\label{cor:Hong}
Let $G$ be a connected graph with $n\ge2$ vertices, $m$ edges, and clique number $\omega$. Then
\[
\rNF(G)\le\demax^{2}\min\Bigl\{\sqrt{2m-n+1},\ \sqrt{2m\bigl(1-\tfrac1\omega\bigr)}\Bigr\}.
\]
More precisely:
\begin{enumerate}[label=(\roman*)]
\item $\rNF(G)\le\demax^{2}\sqrt{2m-n+1}$, with equality if and only if $G\cong K_n$ or
$G\cong K_{1,n-1}$;
\item $\rNF(G)\le\demax^{2}\sqrt{2m(1-1/\omega)}$, with equality if and only if $G$ is a complete
bipartite graph or a complete regular $\omega$-partite graph with $\omega\ge3$.
\end{enumerate}
Consequently, equality holds in the combined bound if and only if $G$ is a complete bipartite graph
or a complete regular multipartite graph.
\end{corollary}

\begin{proof}
(i) By Hong's inequality, $\rho(G)\le\sqrt{2m-n+1}$ for every connected graph, with equality if
and only if $G\cong K_n$ or $G\cong K_{1,n-1}$ \cite{Hong88}; combining this with
Theorem~\ref{thm:delta-rho} gives the bound. Equality requires equality in both estimates, that is,
$\delta$ constant together with $G\cong K_n$ or $G\cong K_{1,n-1}$. Both graphs satisfy the former
condition: $\delta\equiv(n-1)^{2}$ for $K_n$ and $\delta\equiv n-1$ for $K_{1,n-1}$. Conversely,
$\rNF(K_n)=(n-1)^{5}=(n-1)^{4}\sqrt{n^{2}-2n+1}$ and $\rNF(K_{1,n-1})=(n-1)^{2}\sqrt{n-1}$ by
Proposition~\ref{prop:families}, and in both cases the stated bound takes the same value.

(ii) By Nikiforov's inequality \cite{Nik02}, $\rho(G)\le\sqrt{2m(1-1/\omega)}$ for every graph
with $m$ edges and clique number $\omega$; for a graph without isolated vertices, equality holds if
and only if either $\omega=2$ and $G$ is complete bipartite, or $\omega\ge3$ and $G$ is a complete
regular $\omega$-partite graph \cite[Theorem~2]{Nik09}. Combining this with
Theorem~\ref{thm:delta-rho} gives the bound. Equality requires equality in both estimates, that is,
$\delta$ constant together with one of the two structures just described. Both kinds of graphs have
constant $\delta$: $K_{p,q}$ is bipartite semiregular with $\delta\equiv pq$, and a complete regular
$\omega$-partite graph is regular. Conversely, for $K_{p,q}$ we have $\omega=2$, $m=pq$,
$\demax=pq$ and $\rNF(K_{p,q})=(pq)^{2}\sqrt{pq}=\demax^{2}\sqrt{m}$; for $K_{t,\dots,t}$ with
$\omega\ge3$ parts the graph is $r$-regular with $r=(\omega-1)t$, so $\demax=r^{2}$,
$2m(1-1/\omega)=\omega t\cdot r\cdot\frac{\omega-1}{\omega}=r^{2}$ and
$\rNF(G)=r^{5}=\demax^{2}\sqrt{2m(1-1/\omega)}$ by Proposition~\ref{prop:families}(iv).

Finally, $\rNF(G)$ is at most each of the two bounds, so it equals their minimum if and only if it
equals one of them, i.e.\ if and only if equality holds in (i) or in (ii). Since $K_2$ and
$K_{1,n-1}$ are complete bipartite and $K_n$ ($n\ge3$) is complete regular $n$-partite, the graphs
in (i) already occur in (ii), and the combined equality case is that of (ii).
\end{proof}

\begin{remark}\label{rem:Hong-Nik}
The two estimates of Corollary~\ref{cor:Hong} are not comparable. If $G$ is bipartite, then
$\omega=2$ and $2m(1-1/\omega)=m\le2m-n+1$ because $m\ge n-1$, so (ii) is never weaker than (i),
with equality exactly for trees; in this case (ii) also refines Theorem~\ref{thm:bipartite}, since
$\demax\le m$ by Lemma~\ref{lem:bip-delta}. For instance, for $C_6$ part (ii) gives
$16\sqrt6=39.19$ and part (i) gives $16\sqrt7=42.33$, while $\rNF(C_6)=32$; for $K_{2,2,2}$ part
(ii) gives the exact value $1024$, whereas part (i) gives $1115.88$. In the other direction, for
the triangle with one pendant edge (the graph $\Sp4$ of Section~\ref{sec:unicyclic}) one has
$n=m=4$, $\omega=3$ and $2m-n+1=5<\frac{16}{3}=2m(1-1/\omega)$, so (i) is the better estimate:
$25\sqrt5=55.90$ against $57.74$, with $\rNF(\Sp4)=51.51$.
\end{remark}

\begin{remark}\label{rem:Ups-Hong}
The bound of Corollary~\ref{cor:Hong}(i) is \cite[Corollary~3]{BV25}, where equality is asserted
only for $K_n$; as shown above, the star attains it as well. The analogue of
Corollary~\ref{cor:Hong}(i) for the matrix $\Upsilon(G)$ is \cite[Corollary~1]{MBDP22}, where the
upper bound
$2\demax\sqrt{2m-n+1}$ is likewise asserted to be attained only by $K_n$. The star also attains it:
for $G=K_{1,n-1}$ one has $\delta\equiv n-1$, hence
$\rUps(G)=2(n-1)\rho(K_{1,n-1})=2(n-1)\sqrt{n-1}=2\demax\sqrt{2m-n+1}$. The extremal graphs there
are therefore $K_n$ and $K_{1,n-1}$, in accordance with the equality condition in Hong's inequality
as quoted in \cite[Lemma~4]{MBDP22}.
\end{remark}

The estimate of Theorem~\ref{thm:delta-rho} can be improved by observing that only the weights
carried by edges matter. This refinement is the key tool in Section~\ref{sec:unicyclic}.

\begin{lemma}\label{lem:Ww}
Let $G$ be a connected graph and set
\[
W(G)=\max_{\{i,j\}\in E(G)}\delta(i)\delta(j),\qquad w(G)=\min_{\{i,j\}\in E(G)}\delta(i)\delta(j).
\]
Then
\[
W(G)\le\rNF(G)\qquad\text{and}\qquad w(G)\,\rho(G)\le\rNF(G)\le W(G)\,\rho(G),
\]
and either of the last two inequalities is an equality if and only if $G$ is regular or bipartite
semiregular, in which case $w(G)=W(G)=\demax^{2}$ and both are equalities.
\end{lemma}

\begin{proof}
The $(i,j)$-th entry of $\NF(G)$ equals $\delta(i)\delta(j)$ whenever $i\sim j$ and $0$ otherwise,
so, entrywise,
\[
w(G)A(G)\le\NF(G)\le W(G)A(G).
\]
As $G$ is connected, $A(G)$ is irreducible, and Lemma~\ref{lem:PF}(a) gives
$w(G)\rho(G)\le\rNF(G)\le W(G)\rho(G)$, with equality on the right if and only if
$\NF(G)=W(G)A(G)$ and on the left if and only if $\NF(G)=w(G)A(G)$; in either case
$\delta(i)\delta(j)$ is constant on $E(G)$, so Lemma~\ref{lem:const-edge} applies. Conversely, if
$\delta\equiv c$ then $\NF(G)=c^{2}A(G)$ and both inequalities are equalities.

For the first assertion, choose an edge $\{i,j\}$ with $\delta(i)\delta(j)=W(G)$. The principal
submatrix of $\NF(G)$ indexed by $\{i,j\}$ is
$\left[\begin{smallmatrix}0&W(G)\\W(G)&0\end{smallmatrix}\right]$, whose largest eigenvalue is
$W(G)$; now apply Lemma~\ref{lem:interlacing}.
\end{proof}

\begin{remark}\label{rem:Ww}
Since $w(G)\ge\demin^{2}$ and $W(G)\le\demax^{2}$, Lemma~\ref{lem:Ww} refines both bounds of
Theorem~\ref{thm:delta-rho}, and combined with Lemma~\ref{lem:FMS} it refines
Corollary~\ref{cor:delta52}. One has $W(G)=\demax^{2}$ if and only if some two adjacent vertices
both attain $\demax$, and $w(G)=\demin^{2}$ if and only if some two adjacent vertices both attain
$\demin$; in all other cases the refinement is strict. For instance, for the graph $\Sp6$ of
Section~\ref{sec:unicyclic} one has $w=35>\demin^{2}=25$ and $W=49$, and
$49\le\rNF(\Sp6)=110.83\le W\rho(G)=123.19$, while $\demin^{2}\rho(G)=62.85<w\rho(G)=87.99$.
\end{remark}

\begin{remark}\label{rem:rho5-compare}
By Lemma~\ref{lem:Hofmeister} and \eqref{eq:M1-identity} below, $\rho(G)^{2}\ge M_1(G)/n\ge\demin$,
so Theorem~\ref{thm:rho5} gives $\rNF(G)\ge\rho(G)^{5}\ge\demin^{2}\rho(G)$; thus it is never weaker than the lower bound of Theorem~\ref{thm:delta-rho}. It is, however, not comparable with the
lower bound $w(G)\rho(G)$ of Lemma~\ref{lem:Ww}. For $\Sp6$ one has
$\rho^{5}=100.45>w\rho=87.99$. On the other hand, let $D(k,k)$ be the double star obtained by
joining the centres of two copies of $K_{1,k}$ by an edge; for $D(5,5)$ one has $w=66$,
$\rho=(1+\sqrt{21})/2$, $\rho^{5}=169.44<w\rho=184.22$, and $\rNF(D(5,5))=220$. Among the $995$
connected graphs on $2$ to $7$ vertices, $\rho^{5}>w\rho$ for $969$ graphs, $\rho^{5}=w\rho$ for
the $24$ regular or bipartite semiregular ones, and $\rho^{5}<w\rho$ for the remaining $2$.
\end{remark}

\section{Bipartite graphs and trees}\label{sec:bipartite}

\begin{lemma}\label{lem:bip-delta}
Let $G$ be a connected bipartite graph with bipartition $(X,Y)$ and $m$ edges. Then $\demax\le m$.
Moreover, for $v\in V(G)$ we have $\delta(v)=m$ if and only if $v$ is adjacent to every vertex of the
part not containing $v$.
\end{lemma}

\begin{proof}
Since $G$ is bipartite, $\sum_{i\in X}d(i)=\sum_{j\in Y}d(j)=m$. If $v\in X$, then all neighbours
of $v$ lie in $Y$, whence $\delta(v)\le\sum_{j\in Y}d(j)=m$, with equality if and only if $N_G(v)=Y$. The other case is similar.
\end{proof}

\begin{theorem}\label{thm:bipartite}
Let $G$ be a connected bipartite graph with $m$ edges. Then
\[
\rNF(G)\le m^{2}\sqrt{m},
\]
and equality holds if and only if $G\cong K_{p,q}$ with $pq=m$.
\end{theorem}

\begin{proof}
By Corollary~\ref{cor:delta52} and Lemma~\ref{lem:bip-delta},
$\rNF(G)\le\demax^{5/2}\le m^{2}\sqrt{m}$. If equality holds, then equality holds in
Corollary~\ref{cor:delta52}, so $\delta\equiv\demax$, and moreover $\demax=m$. Consequently,
$\delta(v)=m$ for every $v\in V(G)$, and by Lemma~\ref{lem:bip-delta} every vertex is adjacent to all
vertices of the opposite part. Hence $G\cong K_{p,q}$ with $p=|X|$, $q=|Y|$ and $pq=m$. Conversely,
$\rNF(K_{p,q})=(pq)^{2}\sqrt{pq}=m^{2}\sqrt{m}$ by Proposition~\ref{prop:families}(ii).
\end{proof}

To obtain a lower bound for trees, we need the following notation and two lemmas. For integers $a,b,c\ge1$, let $\mathcal{T}(a,b,c)$ denote the tree with a vertex $v$ such that $\mathcal{T}(a,b,c)-v$ consists of three paths on $a$, $b$ and $c$ vertices, each having an end vertex adjacent to $v$.
Let $H_7$ denote the tree obtained from the path $P_3$ by attaching two pendant vertices to each
of its two end vertices, and let $D(2,2)$ be the double star of Remark~\ref{rem:rho5-compare}.
For a subgraph $H$ of a graph $G$ and $i\in V(H)$, let $d_H(i)$, $\delta_H(i)$ and $d_G(i)$,
$\delta_G(i)$ denote the degree and the neighbourhood degree sum of $i$ in $H$ and in $G$.

\begin{lemma}\label{lem:subgraph}
Let $G$ be a connected graph and let $H$ be a connected subgraph of $G$ with at least two
vertices. Then $\rNF(H)\le\rNF(G)$.
\end{lemma}

\begin{proof}
For every $i\in V(H)$ we have $N_H(i)\subseteq N_G(i)$ and $d_H(j)\le d_G(j)$ for all $j$, and so
$\delta_H(i)\le\delta_G(i)$. Hence $\NF(H)\le\NF(G)[V(H)]$ entrywise. Since $H$ is connected, the
matrix $\NF(G)[V(H)]$ is irreducible. By Lemma~\ref{lem:PF}(a) and
Lemma~\ref{lem:interlacing},
\[
\rNF(H)\le\rho\bigl(\NF(G)[V(H)]\bigr)\le\rNF(G).\qedhere
\]
\end{proof}

\begin{lemma}\label{lem:tree-32}
\leavevmode
\begin{enumerate}[label=(\roman*)]
\item $\rNF(P_n)<32$ for every $n\ge2$.
\item Each of the trees $K_{1,4}$, $T(1,2,2)$, $D(2,2)$, $H_7$ and $T(1,1,4)$ has
$\rNF\ge32$.
\end{enumerate}
\end{lemma}

\begin{proof}
(i) Every vertex of $P_n$ has $\delta\le4$, and so $\NF(P_n)\le16A(P_n)$ entrywise. Since
$\rho(P_n)=2\cos\frac{\pi}{n+1}<2$, Lemma~\ref{lem:PF}(a) gives
$\rNF(P_n)\le16\rho(P_n)<32$.

(ii) We have $\rNF(K_{1,4})=32$ by Proposition~\ref{prop:families}(ii). For each of the other
four trees, let $z$ be the positive vector given below. Then $\NF z\ge32z$ entrywise, so
$z^{T}\NF z\ge32\,z^{T}z$, and Lemma~\ref{lem:RR} gives $\rNF\ge32$.
\begin{enumerate}[label=(\alph*)]
\item $\mathcal{T}(1,2,2)$. The vertex of degree $3$ has $\delta=5$, its pendant neighbour has $\delta=3$,
the two vertices of degree $2$ have $\delta=4$, and the other two pendant vertices have
$\delta=2$. Let $z$ take the values $20$, $9$, $13$ and $3$ at these vertices, respectively.
\item $D(2,2)$. The two vertices of degree $3$ have $\delta=5$, and the pendant vertices have
$\delta=3$. Let $z$ take the value $3$ at the vertices of degree $3$ and $1$ at the pendant
vertices.
\item $H_7$. The two vertices of degree $3$ have $\delta=4$, the vertex of degree $2$ has
$\delta=6$, and the pendant vertices have $\delta=3$. Let $z$ take the values $4$, $6$ and $1$ at
these vertices, respectively.
\item $\mathcal{T}(1,1,4)$. Let $v$ be the vertex of degree $3$ and let $vp_1p_2p_3p_4$ be the path of
length $4$ starting at $v$. Then $\delta(v)=4$, the two pendant neighbours of $v$ have
$\delta=3$, and $\delta(p_1),\delta(p_2),\delta(p_3),\delta(p_4)$ are $5,4,3,2$. Let $z$ take the
value $734$ at $v$, the value $275$ at the two pendant neighbours of $v$, and the values $845$,
$618$, $240$ and $45$ at $p_1$, $p_2$, $p_3$ and $p_4$.\qedhere
\end{enumerate}
\end{proof}

\begin{theorem}\label{thm:trees}
Let $G$ be a tree on $n\ge2$ vertices. Then
\[
\rNF(P_n)\le\rNF(G)\le(n-1)^{2}\sqrt{n-1}.
\]
Equality holds in the first inequality if and only if $G\cong P_n$, and in the second inequality
if and only if $G\cong K_{1,n-1}$.
\end{theorem}

\begin{proof}
We first prove the upper bound. A tree on $n$ vertices is bipartite and has $m=n-1$ edges, so
Theorem~\ref{thm:bipartite} gives $\rNF(G)\le(n-1)^{2}\sqrt{n-1}$, with equality if and only if
$G\cong K_{p,q}$ where $pq=n-1$ and $p+q=n$. These two equations force $\{p,q\}=\{1,n-1\}$, i.e.\
$G\cong K_{1,n-1}$, and for the star equality indeed holds by
Proposition~\ref{prop:families}(ii).

Next we prove the lower bound. For $n\le3$ we have $G\cong P_n$. Let $n\ge4$ and suppose that
$G\not\cong P_n$. Then $G$ has a vertex $v$ with $d(v)\ge3$. We show that $\rNF(G)>\rNF(P_n)$,
using Lemma~\ref{lem:subgraph} throughout.

Let $n=4$. Then $G\cong K_{1,3}$, and $\rNF(K_{1,3})=9\sqrt3$ by
Proposition~\ref{prop:families}(ii). The edge weights of $\NF(P_4)$ are $6,9,6$, and for
$z=(1,2,2,1)^{T}$ we have $\NF(P_4)z=12z$. By Lemma~\ref{lem:PF}(c), $\rNF(P_4)\le12<9\sqrt3$.

Let $n\in\{5,6\}$. If $d(v)\ge4$, then $G$ contains $K_{1,4}$, and so $\rNF(G)\ge32$. Otherwise
$d(v)=3$, and since $n\ge5$, some neighbour of $v$ has another neighbour. Hence $G$ contains
$\mathcal{T}(1,1,2)$. In $\mathcal{T}(1,1,2)$, the vertex of degree $3$ has $\delta=4$, its two pendant neighbours
have $\delta=3$, the vertex of degree $2$ has $\delta=4$, and the remaining pendant vertex has
$\delta=2$. The vector $z$ taking the values $4$, $2$, $3$ and $1$ at these vertices satisfies
$\NF(T(1,1,2))z=24z$, and so $\rNF(G)\ge\rNF(\mathcal{T}(1,1,2))\ge24$ by Lemma~\ref{lem:RR}. On the other
hand, the edge weights of $\NF(P_6)$ are $6,12,16,12,6$, and for $z=(1,3,5,5,3,1)^{T}$ we have
$\NF(P_6)z\le23.2\,z$ entrywise. By Lemma~\ref{lem:PF}(c), $\rNF(P_5)\le\rNF(P_6)\le23.2<24$.

Let $n\ge7$. By Lemma~\ref{lem:tree-32}, it suffices to show that $G$ contains one of the trees
$K_{1,4}$, $T(1,2,2)$, $D(2,2)$, $H_7$ and $T(1,1,4)$, since then
$\rNF(G)\ge32>\rNF(P_n)$. If $d(v)\ge4$, then $G$ contains $K_{1,4}$. So let the maximum degree
of $G$ be $3$. If some vertex of degree $3$ has two neighbours that are not pendant vertices,
then $G$ contains $T(1,2,2)$. If two vertices of degree $3$ are adjacent, then $G$ contains
$D(2,2)$. Otherwise, since $n\ge7$, every vertex of degree $3$ is adjacent to exactly two pendant
vertices, and no two vertices of degree $3$ are adjacent. Deleting the pendant vertices adjacent
to vertices of degree $3$ gives a tree of maximum degree at most $2$, that is, a path, and every
vertex of degree $3$ in $G$ becomes an end vertex of this path. Hence $G$ is obtained from a path
by attaching two pendant vertices to one or to both of its end vertices. In the first case $G$
contains $T(1,1,4)$, since $n\ge7$. In the second case $G\cong H_7$ if $n=7$, and $G$ contains
$T(1,1,4)$ if $n\ge8$. This completes the proof.
\end{proof}

\begin{remark}\label{rem:Ups-trees}
Since $m\le\lfloor n^{2}/4\rfloor$ for a bipartite graph on $n$ vertices, Theorem~\ref{thm:bipartite} gives $\rNF(G)\le\lfloor n^{2}/4\rfloor^{5/2}$, with equality if and only if $G\cong K_{\lfloor n/2\rfloor,\lceil n/2\rceil}$. For even $n$ this is \cite[Theorem~8]{BV25}. Theorem~\ref{thm:bipartite} also covers odd $n$ and, being expressed in terms of $m$, is much sharper for sparse graphs: for $K_{1,9}$ it gives the exact value $243$,
whereas \cite[Theorem~8]{BV25} gives $5^{5}=3125$. The analogue of the upper bound in Theorem~\ref{thm:trees} for
$\Upsilon(G)$ is \cite[Lemma~21]{BV24}, deduced there from \cite{ZGJ23}.
\end{remark}

\begin{remark}\label{rem:trees-gap}
The proof shows that $\rNF(G)\ge32$ for every tree $G\not\cong P_n$ on $n\ge7$ vertices, whereas
$\rNF(P_n)<32=\rNF(C_n)$ by Lemma~\ref{lem:tree-32}(i). The bound in
Lemma~\ref{lem:tree-32}(i) is best possible. Indeed, let the vertices of $P_n$ be $1,2,\ldots,n$ in order along the path. For $n\ge6$ the principal submatrix of $\NF(P_n)$ indexed by the vertices $3,4,\ldots,n-2$ is $16A(P_{n-4})$, and so Lemma~\ref{lem:interlacing} gives $\rNF(P_n)\ge32\cos\frac{\pi}{n-3}$,
which tends to $32$.
\end{remark}
\section{Unicyclic graphs}\label{sec:unicyclic}

Let $\mathcal{U}_n$ denote the set of connected unicyclic graphs on $n\ge3$ vertices. Let $\Sp n$
be the graph obtained from the star $K_{1,n-1}$ by joining two of its pendant vertices by an edge.
Equivalently, $\Sp n$ is obtained from the triangle $C_3$ by attaching $n-3$ pendant vertices to
one of its vertices. In particular, $\Sp 3=C_3$. In this section, we determine the minimum and the
maximum of $\rNF(G)$ over $\mathcal{U}_n$ and characterise the graphs attaining them. For the adjacency spectral radius, $C_n$ is the unique minimiser over $\mathcal{U}_n$
and $S_n^+$ is the unique maximiser~\cite{BS86, Sim87}; see also~\cite{Guo05}.
Theorem~\ref{thm:unicyclic} below is the analogue for $N_F(G)$. We first
prove two lemmas.

\begin{lemma}\label{lem:unicyclic-relocation}
Let $G$ be a connected graph and let $v\in V(G)$ with $N_G(v)=\{u,w_1,\ldots,w_q\}$, where
$q\ge1$, each $w_i$ is a pendant vertex and $d(u)\ge2$. Let $G'$ be obtained from $G$ by replacing
the edges $\{v,w_i\}$ with $\{u,w_i\}$, $1\le i\le q$. Then $\rNF(G')>\rNF(G)$.
\end{lemma}

\begin{proof}
Let $S=N_G(u)\setminus\{v\}$, $p=|S|\ge1$, and let $\delta'$ denote the neighbourhood degree sum
in $G'$. Since only the degrees of $u$ and $v$ change, we have
\[
\delta'(u)=\delta(u),\qquad \delta'(v)=\delta(v)=p+q+1,\qquad \delta'(w_i)=\delta(w_i)+p,
\]
and $\delta'(t)\ge\delta(t)$ for every other vertex $t$. Moreover,
$\delta(u)-\delta(v)=\sum_{s\in S}\bigl(d(s)-1\bigr)\ge0$.

Let $X$ be the unit Perron vector of $\NF(G)$, and let $Y$ be obtained from $X$ by setting
$y_u=\max\{x_u,x_v\}$ and $y_v=\min\{x_u,x_v\}$. Then $Y^{T}Y=1$ and $y_uy_v=x_ux_v$, so the edge
$\{u,v\}$ contributes equally to both quadratic forms below. By Lemma~\ref{lem:RR},
\[
\rNF(G')-\rNF(G)\ge Y^{T}\NF(G')Y-X^{T}\NF(G)X
=2\sum_{i=1}^{q}x_{w_i}\bigl(\delta(u)\delta'(w_i)\,y_u-\delta(v)\delta(w_i)\,x_v\bigr)+2\Sigma,
\]
where $\Sigma=\sum\bigl(\delta'(s)\delta'(t)\,y_sy_t-\delta(s)\delta(t)\,x_sx_t\bigr)$ and the sum
is taken over the edges $\{s,t\}\ne\{u,v\}$ common to $G$ and $G'$. Every term of $\Sigma$ is
nonnegative, since $\delta'\ge\delta$ at the end vertices of these edges and $y_u\ge x_u$. Every term of the
first sum is positive, since $\delta(u)\ge\delta(v)$, $\delta'(w_i)>\delta(w_i)$ and
$y_u\ge x_v>0$. Hence $\rNF(G')>\rNF(G)$.
\end{proof}

\begin{lemma}\label{lem:unicyclic-weighted}
Let $n\ge4$, $r>0$ and $3\le g\le n$, and let $H$ be obtained from the cycle $C_g$ by attaching
$n-g$ pendant vertices to vertices of $C_g$. Let $B_r(H)$ be the matrix obtained from $A(H)$ by
replacing the entries corresponding to the pendant edges by $r$. Then
\[
\rho(B_r(H))\le\rho(B_r(\Sp n)),
\]
with equality if and only if $H\cong\Sp n$.
\end{lemma}

\begin{proof}
For $t>1$, define
\[
h(t)=t^{2}-\frac{2t}{t-1}=t^{2}-2-\frac{2}{t-1}.
\]
Since $h'(t)=2t+2/(t-1)^{2}>0$, the function $h$ is strictly increasing on $(1,\infty)$. Let
$s=r^{2}$ and $\Lambda=\rho(B_r(\Sp n))$.
Since $A(C_3)$ is a principal submatrix of $B_r(\Sp n)$, Lemma~\ref{lem:interlacing} gives
$\Lambda\ge2$. Let $z$ be the entry of the Perron vector of $B_r(\Sp n)$ corresponding to the
vertex of degree $n-1$. By symmetry, the entries corresponding to the other two vertices of the triangle are
$z/(\Lambda-1)$, and the entry corresponding to each pendant vertex is $rz/\Lambda$. The
eigenvalue equation at the vertex of degree $n-1$ then gives
\begin{equation}\label{eq:wt-S}
h(\Lambda)=(n-3)s.
\end{equation}

Now let $\lambda=\rho(B_r(H))$. Since $A(C_g)$ is a principal submatrix of $B_r(H)$,
Lemma~\ref{lem:interlacing} gives $\lambda\ge2$. Write $C_g=v_1v_2\cdots v_gv_1$, let $t_i$ be the
number of pendant vertices adjacent to $v_i$, so that $\sum_{i=1}^{g}t_i=n-g$, and let $x_i$ be
the entry of the Perron vector of $B_r(H)$ corresponding to $v_i$. Without loss of generality, we
assume that $x_1=\max_ix_i$. The entry corresponding to a pendant vertex adjacent to $v_i$ is
$rx_i/\lambda$, and so the eigenvalue equation at $v_i$ is
\begin{equation}\label{eq:wt-cycle}
\lambda x_i=x_{i-1}+x_{i+1}+\frac{st_i}{\lambda}\,x_i\qquad(1\le i\le g),
\end{equation}
with indices taken modulo $g$. Let $S=x_1+\cdots+x_g$. Summing \eqref{eq:wt-cycle} over all $i$,
we obtain
\begin{equation}\label{eq:wt-sum}
(\lambda-2)S=\frac{s}{\lambda}\sum_{i=1}^{g}t_ix_i\le\frac{(n-g)s}{\lambda}\,x_1.
\end{equation}
Summing \eqref{eq:wt-cycle} over $i=2,\ldots,g$ and using $x_2+x_g\le S-x_1$, we obtain
\begin{equation}\label{eq:wt-part}
\lambda(S-x_1)=2S-x_2-x_g+\frac{s}{\lambda}\sum_{i=2}^{g}t_ix_i\ge S+x_1,
\end{equation}
that is, $(\lambda-1)S\ge(\lambda+1)x_1$. Since $\lambda\ge2$, \eqref{eq:wt-sum} and
\eqref{eq:wt-part} give
\[
\frac{(\lambda-2)(\lambda+1)}{\lambda-1}\,x_1\le(\lambda-2)S\le\frac{(n-g)s}{\lambda}\,x_1.
\]
Multiplying by $\lambda/x_1$ and using
$\lambda(\lambda-2)(\lambda+1)/(\lambda-1)=h(\lambda)$, we obtain
\begin{equation}\label{eq:wt-final}
h(\lambda)\le(n-g)s\le(n-3)s.
\end{equation}
By \eqref{eq:wt-S} we have $h(\lambda)\le h(\Lambda)$, and since $h$ is strictly increasing,
$\lambda\le\Lambda$.

If $\lambda=\Lambda$, then equality holds throughout. Note that $\lambda>2$, because
$h(\lambda)=(n-3)s>0=h(2)$. Equality in \eqref{eq:wt-final} gives $g=3$, and equality in
\eqref{eq:wt-part} gives $t_i=0$ for $i\ge2$. Thus all pendant vertices are adjacent to $v_1$, and
$H\cong\Sp n$. Conversely, if $H\cong\Sp n$, then $\rho(B_r(H))=\rho(B_r(\Sp n))$.
\end{proof}

We can now state the main result of this section.

\begin{theorem}\label{thm:unicyclic}
Let $G$ be a connected unicyclic graph on $n\ge3$ vertices. Then
\[
32=\rNF(C_n)\le\rNF(G)\le\rNF(\Sp n).
\]
Equality holds in the first inequality if and only if $G\cong C_n$, and in the second inequality
if and only if $G\cong\Sp n$.
\end{theorem}

\begin{proof}
We first prove the lower bound. Let $C$ be the unique cycle of $G$ and let $g=|V(C)|$. Every
vertex of $C$ has two neighbours on $C$, each of degree at least $2$. Hence $\delta(i)\ge4$ for every $i\in V(C)$. Since $G$ is
unicyclic, $C$ has no chords, so the edges of $G$ with both end vertices in $V(C)$ are exactly the
$g$ edges of $C$. Let $X=(x_1,\ldots,x_n)^{T}$ with $x_i=1$ if $i\in V(C)$ and $x_i=0$ otherwise.
By Lemma~\ref{lem:RR},
\[
\rNF(G)\ge\frac{X^{T}\NF(G)X}{X^{T}X}=\frac{2}{g}\sum_{\{i,j\}\in E(C)}\delta(i)\delta(j)
\ge\frac{2}{g}\cdot16g=32.
\]
Now suppose that $G\not\cong C_n$. Then $g<n$. Since $G$ is connected, there exists a vertex
$i\in V(C)$ that has a neighbour outside $C$, and so $d(i)\ge3$. As the two neighbours of $i$ on
$C$ have degree at least $2$, we have $\delta(i)\ge5$. Consequently, $\delta(i)\delta(j)\ge20>16$
for each neighbour $j$ of $i$ on $C$, and the last inequality above is strict. Hence
$\rNF(G)>32$. Finally, $\rNF(C_n)=32$ by Proposition~\ref{prop:families}(iii). This proves the
lower bound and its equality case.

Next we prove the upper bound. For $n=3$ we have $G\cong C_3=\Sp3$. Let $n\ge4$. Since
$\mathcal{U}_n$ is finite, we may choose a graph $G_0\in\mathcal{U}_n$ with the maximum value of
$\rNF$. Let $C$ be the unique cycle of $G_0$. We show that every vertex outside $C$ is a pendant
vertex. If not, let $v$ be a non-pendant vertex outside $C$ at maximum
distance from $C$, and let $u$ be the neighbour of $v$ on a shortest path from $v$ to $C$. Then
all other neighbours of $v$ are pendant vertices, and $d(u)\ge2$. By
Lemma~\ref{lem:unicyclic-relocation}, moving these pendant vertices from $v$ to $u$ gives a graph
in $\mathcal{U}_n$ with a larger value of $\rNF$, a contradiction. Hence $G_0$ is obtained from $C$
by attaching pendant vertices to vertices of $C$.

Let $a=n+1$ and $r=(n-1)/(n+1)$. Since $G_0$ has $n$ edges, $\sum_{j\in V(G_0)}d(j)=2n$, and so
every vertex $i$ of $G_0$ satisfies
\[
\delta(i)=n+1-\sum_{j\notin N_{G_0}[i]}\bigl(d(j)-1\bigr)\le a.
\]
If $i$ is a pendant vertex with neighbour $k$, then $\delta(i)=d(k)\le n-1=ar$. Hence
$\delta(i)\delta(j)\le a^{2}$ on the edges of $C$ and $\delta(i)\delta(j)\le a^{2}r$ on the pendant
edges, that is, $\NF(G_0)\le a^{2}B_r(G_0)$ entrywise. In $\Sp n$, the three vertices of the
triangle have $\delta=n+1$ and the pendant vertices have $\delta=n-1$, and so $\NF(\Sp n)=a^{2}B_r(\Sp n)$.
By Lemma~\ref{lem:PF}(a) and Lemma~\ref{lem:unicyclic-weighted},
\[
\rNF(G_0)\le a^{2}\rho(B_r(G_0))\le a^{2}\rho(B_r(\Sp n))=\rNF(\Sp n).
\]
By the choice of $G_0$, equality holds throughout, and Lemma~\ref{lem:unicyclic-weighted} gives
$G_0\cong\Sp n$. Thus $\Sp n$ is the unique graph in $\mathcal{U}_n$ with the maximum value of
$\rNF$. Hence $\rNF(G)\le\rNF(\Sp n)$, with equality if and only if $G\cong\Sp n$. This completes
the proof.
\end{proof}

\begin{remark}\label{rem:Pstar}
For the neighbourhood first Zagreb matrix $\Upsilon(G)$, the analogue of
Theorem~\ref{thm:unicyclic} was proved in
\cite[Lemma~22]{BV24} for $n\ge7$ by using the results of \cite{ZGJ23} on degree-based weighted
adjacency matrices. These results do not apply to $\NF(G)$, since the weight $\delta(i)\delta(j)$
is not a function of $d(i)$ and $d(j)$ alone. Our proof uses the weights
$\delta(i)\delta(j)$ directly and holds for all $n\ge3$.
\end{remark}



\section{Bounds involving the first Zagreb index}\label{sec:M1}

Recall that the first Zagreb index of $G$ is $M_1(G)=\sum_{i=1}^{n}d(i)^{2}$ \cite{GT72}. Counting
the pairs (edge, incident vertex) in two ways gives the identity
\begin{equation}\label{eq:M1-identity}
\sum_{i=1}^{n}\delta(i)=\sum_{i=1}^{n}\sum_{j\sim i}d(j)=\sum_{j=1}^{n}d(j)^{2}=M_1(G).
\end{equation}

\begin{lemma}\label{lem:delta-M1}
Let $G$ be a connected graph on $n\ge2$ vertices. Then
\[
\demax\le\tfrac12M_1(G),
\]
with equality if and only if $G\cong K_2$.
\end{lemma}

\begin{proof}
Choose $v\in V(G)$ with $\delta(v)=\demax$. For every $j\sim v$ we have
\[
\delta(j)=\sum_{k\sim j}d(k)\ge d(v)+\bigl(d(j)-1\bigr),
\]
because the neighbour $v$ contributes $d(v)$ and each of the remaining $d(j)-1$ neighbours
contributes at least $1$. Summing over $j\sim v$,
\[
\sum_{j\sim v}\delta(j)\ge d(v)^{2}+\sum_{j\sim v}d(j)-d(v)=d(v)^{2}-d(v)+\delta(v).
\]
Hence, using \eqref{eq:M1-identity} and $\delta(i)\ge1$ for all $i$,
\[
M_1(G)=\sum_{i=1}^{n}\delta(i)\ge\delta(v)+\sum_{j\sim v}\delta(j)\ge2\delta(v)+d(v)\bigl(d(v)-1\bigr)\ge2\demax.
\]
If equality holds, then $d(v)(d(v)-1)=0$, so $d(v)=1$; moreover no vertex lies outside $N_G[v]$,
since every vertex $i$ has $\delta(i)\ge1$. Thus $V(G)=\{v,j\}$ where $j$ is the unique neighbour
of $v$, i.e.\ $G\cong K_2$. Conversely, for $K_2$ we have $\demax=1$ and $M_1=2$.
\end{proof}

For the lower bound we use the following classical estimate of Hofmeister \cite{Hof88}, whose short
proof we include for completeness.

\begin{lemma}[\cite{Hof88}]\label{lem:Hofmeister}
Let $G$ be a connected graph on $n$ vertices. Then $\rho(G)\ge\sqrt{M_1(G)/n}$, with equality if and
only if $G$ is regular or bipartite semiregular.
\end{lemma}

\begin{proof}
Let $\one$ be the all-ones vector. Since
$\one^{T}A(G)^{2}\one=\sum_i\sum_{j\sim i}d(j)=\sum_i\delta(i)=M_1(G)$ by \eqref{eq:M1-identity},
the Rayleigh--Ritz inequality applied to $A(G)^{2}$ gives
\[
\rho(G)^{2}=\lambda_{\max}\bigl(A(G)^{2}\bigr)\ge\frac{\one^{T}A(G)^{2}\one}{\one^{T}\one}=\frac{M_1(G)}{n}.
\]
Equality holds if and only if $\one$ is an eigenvector of $A(G)^{2}$, that is, if and only if
$\delta$ is constant, which by Lemma~\ref{lem:delta-const} happens exactly when $G$ is regular or
bipartite semiregular.
\end{proof}

\begin{theorem}\label{thm:M1}
Let $G$ be a connected graph on $n\ge2$ vertices. Then
\[
\rNF(G)\ge w(G)\sqrt{\frac{M_1(G)}{n}}\ge\demin^{2}\sqrt{\frac{M_1(G)}{n}},
\]
and the first inequality is an equality if and only if $G$ is regular or bipartite semiregular, in
which case both are equalities.
\end{theorem}

\begin{proof}
Combine Lemma~\ref{lem:Ww} with Lemma~\ref{lem:Hofmeister}; the second inequality holds because
$w(G)\ge\demin^{2}$. Equality in the first estimate requires equality in both ingredients, and each
of them holds precisely when $\delta$ is constant (Lemmas~\ref{lem:const-edge}
and~\ref{lem:delta-const}); conversely, if $\delta\equiv c$ then
$\rNF(G)=c^{2}\rho(G)=c^{2}\sqrt{M_1(G)/n}$ by Proposition~\ref{prop:factor} and
Lemma~\ref{lem:Hofmeister}. Finally $w(G)=\demin^{2}$ exactly when some two adjacent vertices both
attain $\demin$, which is the case when $\delta$ is constant.
\end{proof}

The weaker estimate $\rNF(G)\ge\demin^{2}\sqrt{M_1(G)/n}$ is the lower bound of
\cite[Corollary~4]{BV25}, where equality is asserted only for regular graphs; by
Theorem~\ref{thm:M1} the bipartite semiregular graphs, for instance the stars, attain it as well.
Theorem~\ref{thm:rho5} gives a second lower bound in terms of $M_1$.

\begin{corollary}\label{cor:M1-rho5}
Let $G$ be a connected graph on $n\ge2$ vertices. Then
\[
\rNF(G)\ge\rho(G)^{5}\ge\Bigl(\frac{M_1(G)}{n}\Bigr)^{5/2},
\]
and equality holds in either inequality if and only if $G$ is regular or bipartite semiregular.
\end{corollary}

\begin{proof}
The first inequality and its equality case are Theorem~\ref{thm:rho5}; the second inequality and
its equality case follow from Lemma~\ref{lem:Hofmeister}.
\end{proof}

The bounds of Theorem~\ref{thm:M1} and Corollary~\ref{cor:M1-rho5} are not comparable. For $\Sp6$
one has $(M_1/n)^{5/2}=88.18>w\sqrt{M_1/n}=85.73$, whereas for the double star $D(5,5)$ of
Remark~\ref{rem:rho5-compare} one has $(M_1/n)^{5/2}=122.06<w\sqrt{M_1/n}=172.53$.

\begin{remark}\label{rem:M1-upper}
Lemma~\ref{lem:delta-M1} also yields upper bounds, namely
\[
\rho(G)\le\sqrt{\tfrac12M_1(G)}\qquad\text{and}\qquad\rNF(G)\le\demax^{2}\sqrt{\tfrac12M_1(G)},
\]
each with equality if and only if $G\cong K_2$. We record them because $M_1$ is a standard
molecular descriptor and estimates expressed through it are customary in this literature, but we
stress that they carry no information beyond Section~\ref{sec:general}: since
$\demax\le\frac12M_1(G)$ with equality only for $K_2$, the first is weaker than
Lemma~\ref{lem:FMS} and the second is weaker than Corollary~\ref{cor:delta52} for every graph other
than $K_2$. It is the lower bounds of Theorem~\ref{thm:M1} and Corollary~\ref{cor:M1-rho5} that are
of genuine use, and they are not comparable with the upper estimates of Section~\ref{sec:general}.
\end{remark}

\section{Bounds for the neighbourhood second Zagreb energy}\label{sec:energy}

In analogy with the graph energy of Gutman \cite{Gut78}, the neighbourhood second Zagreb energy of
a graph $G$ was defined in \cite{BV25} as
\[
\E(\NF)=\sum_{i=1}^{n}|\lambda_i|,
\]
where $\lambda_1\ge\dots\ge\lambda_n$ are the eigenvalues of $\NF(G)$. Throughout this section $G$
has no isolated vertices, $T(G)$ is as in \eqref{eq:T} and
\[
P(G)=\bigl|\det\NF(G)\bigr|^{2/n}=\Bigl(\prod_{i=1}^{n}\delta(i)\Bigr)^{4/n}\bigl|\det A(G)\bigr|^{2/n},
\]
the second expression being furnished by Proposition~\ref{prop:congruent}. We first record the
estimate of McClelland type \cite{McC71}, which follows from Lemma~\ref{lem:moments} and the
Cauchy--Schwarz inequality:
\begin{equation}\label{eq:McClelland}
\E(\NF)\le\sqrt{2nT(G)}.
\end{equation}

\begin{theorem}[{\cite[Lemma~9]{BV25}}]\label{thm:Kober}
Let $G$ be a graph on $n$ vertices with no isolated vertices. Then
\[
\sqrt{2T(G)+n(n-1)P(G)}\le\E(\NF)\le\sqrt{2(n-1)T(G)+nP(G)}.
\]
Moreover, the upper bound never exceeds \eqref{eq:McClelland}, and it is strictly smaller unless
every $|\lambda_i|$ takes the same value, which for a graph without isolated vertices happens
exactly when $n$ is even and $G\cong\frac n2K_2$.
\end{theorem}

\begin{proof}
For nonnegative real numbers $a_1,\dots,a_n$, Kober's inequality \cite{Kob58} states
\[
\sum_{i=1}^{n}a_i+n(n-1)\Bigl(\prod_{i=1}^{n}a_i\Bigr)^{1/n}\le\Bigl(\sum_{i=1}^{n}\sqrt{a_i}\Bigr)^{2}
\le(n-1)\sum_{i=1}^{n}a_i+n\Bigl(\prod_{i=1}^{n}a_i\Bigr)^{1/n}.
\]
Take $a_i=\lambda_i^{2}$. Then $\sum_ia_i=2T(G)$ by Lemma~\ref{lem:moments},
$\sum_i\sqrt{a_i}=\E(\NF)$ and
$\bigl(\prod_ia_i\bigr)^{1/n}=\bigl(\prod_i\lambda_i^{2}\bigr)^{1/n}=|\det\NF(G)|^{2/n}=P(G)$, which
gives both inequalities.

By the arithmetic--geometric mean inequality, $P(G)\le\frac1n\sum_i\lambda_i^{2}=\frac{2T(G)}{n}$,
whence $2(n-1)T(G)+nP(G)\le2nT(G)$, with equality if and only if all $\lambda_i^{2}$ are equal. In
that case $\NF(G)^{2}=c^{2}I_n$ for some $c>0$, so for $i\ne j$ we get
$0=(\NF^{2})_{ij}=\sum_kn_{ik}n_{kj}$ with all summands nonnegative, and therefore $n_{ik}n_{kj}=0$
for every $k$; thus no two distinct vertices have a common neighbour, so every component of $G$ has
at most two vertices. As $G$ has no isolated vertices, $G\cong\frac n2K_2$, and conversely for this
graph $\delta\equiv1$, $\NF(G)=A(G)$ and all $|\lambda_i|=1$.
\end{proof}

The next bound is of Koolen--Moulton type \cite{KM01}: it exploits any available lower estimate for
$\rNF(G)$, and the sharper the estimate, the better the bound.

\begin{theorem}[Koolen--Moulton type, cf.\ {\cite{KM01}}]\label{thm:KM}
Let $G$ be a graph on $n$ vertices with no isolated vertices, and let $\rho_0$ satisfy
$\sqrt{2T(G)/n}\le\rho_0\le\rNF(G)$. Then
\[
\E(\NF)\le\rho_0+\sqrt{(n-1)\bigl(2T(G)-\rho_0^{2}\bigr)}.
\]
In particular, the choice
\[
\rho_0=\max\Bigl\{\tfrac2nM_2'(G),\ w(G)\rho(G),\ W(G),\ \rho(G)^{5},\ \sqrt{2T(G)/n}\Bigr\}
\]
is admissible when $G$ is connected, and the resulting bound never exceeds \eqref{eq:McClelland}.
\end{theorem}

\begin{proof}
By the Cauchy--Schwarz inequality and Lemma~\ref{lem:moments},
\[
\E(\NF)=\rNF(G)+\sum_{i=2}^{n}|\lambda_i|\le\rNF(G)+\sqrt{(n-1)\sum_{i=2}^{n}\lambda_i^{2}}=g\bigl(\rNF(G)\bigr),
\]
where $g(x)=x+\sqrt{(n-1)(2T(G)-x^{2})}$. On the interval
$\bigl[\sqrt{2T(G)/n},\sqrt{2T(G)}\bigr]$ one has
$g'(x)=1-\frac{(n-1)x}{\sqrt{(n-1)(2T(G)-x^{2})}}\le0$, so $g$ is nonincreasing there. Since
$2T(G)=\sum_i\lambda_i^{2}\le n\rNF(G)^{2}$ we have $\rNF(G)\ge\sqrt{2T(G)/n}$, and
$\rNF(G)\le\sqrt{2T(G)}$ trivially; hence $\rho_0\le\rNF(G)$ gives $g(\rNF(G))\le g(\rho_0)$, which
is the asserted bound.

The five quantities listed are lower bounds for $\rNF(G)$: the first because the Rayleigh quotient
at the all-ones vector equals $\frac1n\one^{T}\NF(G)\one=\frac2nM_2'(G)$, the second and the third
by Lemma~\ref{lem:Ww}, the fourth by Theorem~\ref{thm:rho5}, and the fifth as shown above. Finally
$g\bigl(\sqrt{2T(G)/n}\bigr)=\sqrt{2T(G)/n}+\sqrt{(n-1)^{2}2T(G)/n}=\sqrt{2nT(G)}$, so the bound is
never weaker than \eqref{eq:McClelland}.
\end{proof}

\begin{remark}\label{rem:KM-Kn}
The bound of Theorem~\ref{thm:KM} is attained by $K_n$: there $\delta\equiv(n-1)^{2}$,
$2T=n(n-1)^{9}$ and $\rNF=(n-1)^{5}=\rho_0$, so the right-hand side equals
$2(n-1)^{5}=\E(\NF(K_n))$.
\end{remark}

We now turn to lower bounds, where Proposition~\ref{prop:congruent} yields complete
characterisations of the extremal graphs.

\begin{theorem}\label{thm:energy-lower}
Let $G$ be a graph on $n$ vertices with no isolated vertices. Then
\[
\E(\NF)\ge2\rNF(G),
\]
and $\E(\NF)=2\rNF(G)$ if and only if $G$ is a complete multipartite graph; in particular $G$ is
then connected, since a graph with two components containing edges has at least two positive
$\NF$-eigenvalues. If $G$ is connected, then moreover
\[
2\rNF(G)\ge2\max\bigl\{W(G),\ w(G)\rho(G),\ \rho(G)^{5}\bigr\}.
\]
\end{theorem}

\begin{proof}
By Lemma~\ref{lem:moments} the eigenvalues of $\NF(G)$ sum to zero, so the positive ones and the
negative ones have the same absolute total, giving
$\E(\NF)=2\sum_{\lambda_i>0}\lambda_i\ge2\lambda_1=2\rNF(G)$, with equality if and only if
$\NF(G)$ has exactly one positive eigenvalue. By Proposition~\ref{prop:congruent} this happens if
and only if $A(G)$ has exactly one positive eigenvalue, which by a theorem of Smith \cite{Smi70}
characterises the complete multipartite graphs among graphs without isolated vertices; equivalently,
by Corollary~\ref{cor:inertia}(iii). The last inequality is Lemma~\ref{lem:Ww} together with
Theorem~\ref{thm:rho5}.
\end{proof}

\begin{corollary}\label{cor:unicyclic-energy}
Let $G$ be a connected unicyclic graph. Then $\E(\NF)\ge64$, with equality if and only if
$G\cong C_3$ or $G\cong C_4$.
\end{corollary}

\begin{proof}
By Theorem~\ref{thm:unicyclic}, $\rNF(G)\ge32$, so $\E(\NF)\ge64$ by
Theorem~\ref{thm:energy-lower}. Equality requires $\rNF(G)=32$, hence $G\cong C_n$ by
Theorem~\ref{thm:unicyclic}, and simultaneously that $G$ be complete multipartite; among the cycles
only $C_3=K_{1,1,1}$ and $C_4=K_{2,2}$ are complete multipartite. For both, $\NF(G)=16A(G)$ and
$\E(\NF)=16\,\E(G)=64$.
\end{proof}

\begin{theorem}\label{thm:energy-T}
Let $G$ be a connected graph on $n$ vertices. Then
\[
\E(\NF)\ge\frac{2T(G)}{\rNF(G)}\ge\frac{2T(G)}{W(G)\rho(G)}\ge\frac{2T(G)}{\demax^{5/2}},
\]
and the first inequality is an equality if and only if $G$ is a complete bipartite graph.
\end{theorem}

\begin{proof}
By Lemma~\ref{lem:moments},
\[
2T(G)=\sum_{i=1}^{n}\lambda_i^{2}=\sum_{i=1}^{n}|\lambda_i|\,|\lambda_i|\le\rNF(G)\sum_{i=1}^{n}|\lambda_i|=\rNF(G)\,\E(\NF),
\]
which gives the first inequality. The second follows from Lemma~\ref{lem:Ww}, and the third from
$W(G)\le\demax^{2}$ together with $\rho(G)\le\sqrt{\demax}$ (Lemma~\ref{lem:FMS}), which give
$W(G)\rho(G)\le\demax^{5/2}$.

Equality in the displayed estimate holds if and only if $|\lambda_i|\bigl(\rNF(G)-|\lambda_i|\bigr)=0$
for every $i$, that is, if and only if every eigenvalue of $\NF(G)$ lies in
$\{\rNF(G),0,-\rNF(G)\}$. As $G$ is connected, $\NF(G)$ is irreducible and $\rNF(G)$ is a simple
eigenvalue; since the eigenvalues sum to zero, $-\rNF(G)$ then occurs exactly once as well, so
$\rank\NF(G)=2$. By Proposition~\ref{prop:congruent}, $\rank A(G)=2$, so $G\cong K_{p,q}$ by
Corollary~\ref{cor:inertia}(ii). Conversely, for $G=K_{p,q}$ we have $\NF(G)=(pq)^{2}A(G)$, whence
$\E(\NF)=2(pq)^{5/2}$, $\rNF(G)=(pq)^{5/2}$ and $T(G)=(pq)^{5}$, so equality holds.
\end{proof}

\section*{Acknowledgments}
The second author acknowledges financial support from the ANRF-CRG grant (File Number: CRG/2023/002747).
The authors acknowledge the use of AI-Claude for improving the exposition of the manuscript and refining some of the proofs. All mathematical statements have been independently verified by the authors, who are solely responsible for the correctness of the results presented in this paper. 
\section*{Declarations}

\noindent\textbf{Conflict of interest:} The authors declare that they have no competing interests.



\begin{thebibliography}{99}

\bibitem{BV25}
S. Barik and P. Verma, On spectral properties of neighborhood second Zagreb matrix of graphs,
\emph{Commun. Comb. Optim.} \textbf{10} (2025), no.~2, 275--293.
\url{https://doi.org/10.22049/cco.2023.28424.1546}

\bibitem{BV24}
S. Barik and P. Verma, On the neighborhood first Zagreb spectral radius and the neighborhood first
Zagreb energy of graphs, \emph{AKCE Int. J. Graphs Comb.}, accepted (2024).

\bibitem{BZ01}
A. Berman and X.-D. Zhang, On the spectral radius of graphs with cut vertices,
\emph{J. Combin. Theory Ser. B} \textbf{83} (2001), 233--240.

 \bibitem{BS86}
 R. A. Brualdi and E. S. Solheid, On the spectral radius of complementary acyclic matrices of zeros and ones, \emph{SIAM J. Algebraic Discrete Methods} \textbf{7} (1986), 265--272.

\bibitem{Cao98}
D. S. Cao, Bounds on the eigenvalues and chromatic number, \emph{Linear Algebra Appl.} \textbf{270}
(1998), 1--13.

\bibitem{CL07}
B. Cheng and B. Liu, On the nullity of graphs, \emph{Electron. J. Linear Algebra} \textbf{16}
(2007), 60--67.

\bibitem{CRS10}
D. Cvetkovi\'c, P. Rowlinson and S. Simi\'c, \emph{An Introduction to the Theory of Graph Spectra},
Cambridge University Press, Cambridge, 2010.

\bibitem{DGMMF18}
K. C. Das, I. Gutman, I. Milovanovi\'c, E. Milovanovi\'c and B. Furtula, Degree-based energies of
graphs, \emph{Linear Algebra Appl.} \textbf{554} (2018), 185--204.

\bibitem{FMS93}
O. Favaron, M. Mah\'eo and J.-F. Sacl\'e, Some eigenvalue properties in graphs (conjectures of
Graffiti~II), \emph{Discrete Math.} \textbf{111} (1993), 197--220.

\bibitem{FG15}
B. Furtula and I. Gutman, A forgotten topological index, \emph{J. Math. Chem.} \textbf{53} (2015),
1184--1190.

\bibitem{GH10}
M. Ghorbani and M. A. Hosseinzadeh, A note on Zagreb indices of nanostar dendrimers,
\emph{Optoelectron. Adv. Mater. Rapid Comm.} \textbf{4} (2010), 1877--1880.

\bibitem{Guo05}
S. G. Guo, The spectral radius of unicyclic and bicyclic graphs with $n$ vertices and $k$ pendant
vertices, \emph{Linear Algebra Appl.} \textbf{408} (2005), 78--85.

\bibitem{Gut78}
I. Gutman, The energy of a graph, \emph{Ber. Math.-Statist. Sekt. Forschungsz. Graz} \textbf{103}
(1978), 1--22.

\bibitem{GT72}
I. Gutman and N. Trinajsti\'c, Graph theory and molecular orbitals. Total $\pi$-electron energy of
alternant hydrocarbons, \emph{Chem. Phys. Lett.} \textbf{17} (1972), 535--538.

\bibitem{GRTW75}
I. Gutman, B. Ru\v{s}\v{c}i\'c, N. Trinajsti\'c and C. F. Wilcox, Graph theory and molecular
orbitals. XII. Acyclic polyenes, \emph{J. Chem. Phys.} \textbf{62} (1975), 3399--3405.

\bibitem{Hof88}
M. Hofmeister, Spectral radius and degree sequence, \emph{Math. Nachr.} \textbf{139} (1988),
37--44.

\bibitem{Hong88}
Y. Hong, A bound on the spectral radius of graphs, \emph{Linear Algebra Appl.} \textbf{108} (1988),
135--139.

\bibitem{HJ13}
R. A. Horn and C. R. Johnson, \emph{Matrix Analysis}, 2nd ed., Cambridge University Press,
Cambridge, 2013.

\bibitem{JJG18}
N. Jafari Rad, A. Jahanbani and I. Gutman, Zagreb energy and Zagreb Estrada index of graphs,
\emph{MATCH Commun. Math. Comput. Chem.} \textbf{79} (2018), 371--386.

\bibitem{Kob58}
H. Kober, On the arithmetic and geometric means and on H\"older's inequality, \emph{Proc. Amer.
Math. Soc.} \textbf{9} (1958), 452--459.

\bibitem{KM01}
J. H. Koolen and V. Moulton, Maximal energy graphs, \emph{Adv. Appl. Math.} \textbf{26} (2001),
47--52.

\bibitem{LW21}
X. Li and Z. Wang, Trees with extremal spectral radius of weighted adjacency matrices among trees
weighted by degree-based indices, \emph{Linear Algebra Appl.} \textbf{620} (2021), 61--75.

\bibitem{McC71}
B. J. McClelland, Properties of the latent roots of a matrix: the estimation of $\pi$-electron
energies, \emph{J. Chem. Phys.} \textbf{54} (1971), 640--643.

\bibitem{MBDP22}
S. Mondal, S. Barik, N. De and A. Pal, A note on neighborhood first Zagreb energy and its
significance as a molecular descriptor, \emph{Chemometr. Intell. Lab. Syst.} \textbf{222} (2022),
Article ID 104494.

\bibitem{MDP19}
S. Mondal, N. De and A. Pal, On some new neighborhood degree based indices, \emph{Acta Chem. Iasi}
\textbf{27} (2019), 31--46.

\bibitem{MDP21}
S. Mondal, N. De and A. Pal, On neighborhood Zagreb index of product graphs, \emph{J. Mol. Struct.}
\textbf{1223} (2021), 129210.

\bibitem{Nik02}
V. Nikiforov, Some inequalities for the largest eigenvalue of a graph, \emph{Combin. Probab.
Comput.} \textbf{11} (2002), 179--189.

\bibitem{Nik09}
V. Nikiforov, More spectral bounds on the clique and independence numbers, \emph{J. Combin. Theory
Ser. B} \textbf{99} (2009), 819--826.

\bibitem{Sim87}
S. K. Simi\'c, On the largest eigenvalue of unicyclic graphs, \emph{Publ. Inst. Math. (Beograd)}
\textbf{42} (1987), 13--19.

\bibitem{Smi70}
J. H. Smith, Some properties of the spectrum of a graph, in: \emph{Combinatorial Structures and
Their Applications}, Gordon and Breach, New York, 1970, pp.~403--406.

\bibitem{VOK09}
P. Van Mieghem, J. Omic and R. Kooij, Virus spread in networks, \emph{IEEE/ACM Trans. Netw.}
\textbf{17} (2009), no.~1, 1--14.

\bibitem{ZGJ23}
R. Zheng, X. Guan and X. Jin, Extremal trees and unicyclic graphs with respect to spectral radius
of weighted adjacency matrices with property $P^{*}$, \emph{J. Appl. Math. Comput.} \textbf{69}
(2023), no.~3, 2573--2594.

\end{thebibliography}
\end{document}